\documentclass[a4paper, preprint, review, 11pt]{article}

\usepackage[linesnumbered,ruled,vlined]{algorithm2e}

\SetCommentSty{mycommfont}

\SetKwInput{KwInput}{Input}                
\SetKwInput{KwOutput}{Output}              

\usepackage{caption}
\usepackage{amsmath,amssymb}

\usepackage{adjustbox}
 \usepackage{bbm}
 \usepackage{todonotes}
 \usepackage{graphicx}
\usepackage{color}
\definecolor{cornell-red}{RGB}{179,27,27}
\usepackage{subcaption}
\usepackage{mathtools}
\usepackage{dsfont}
\usepackage{float}
\usepackage{natbib}
    \usepackage{multirow}
 
\usepackage[left=1in,right=1in,top=1in,bottom=1in]{geometry}
\usepackage{rotating}
\usepackage{xcolor}
\usepackage[colorlinks = true,linkcolor = blue,urlcolor  = blue,citecolor = blue,anchorcolor = blue]{hyperref}
\usepackage[flushleft]{threeparttable}

\usepackage{booktabs}
\usepackage{multirow}
\usepackage{rotating,tabularx}

\usepackage{xcolor}

\def\mathcolor#1#{\@mathcolor{#1}}
\def\@mathcolor#1#2#3{%
  \protect\leavevmode
  \begingroup
    \color#1{#2}#3%
  \endgroup
}
\usepackage{pifont}
\usepackage{graphicx}
\usepackage{natbib}
\usepackage{lscape}

 \bibpunct[, ]{(}{)}{,}{a}{}{,}%
 \def\bibfont{\small}%
\begin{document}




\title{Adaptive Large Neighborhood Search Metaheuristic for Vehicle Routing Problem with Multiple Synchronization Constraints and Multiple Trips}


\author{Faisal Alkaabneh$^{1*}$}
\date{\vspace{-5ex}}

\maketitle

\author{$^{1}$Industrial \& Systems Engineering, North Carolina A\&T State University, Greensboro, NC 27401, United States, Email: fmalkaabneh@ncat.edu. $^{*}$\textit{Corresponding Author}}

\abstract{%
This work is motivated by solving a problem faced by big agriculture companies implementing precision agriculture operations for spraying practices using two types of operators, namely a tender tanker and a fleet of sprayers. We model this problem as a vehicle routing problem with multiple synchronization constraints and multiple trips with the objective of minimizing the waiting time of the sprayers and the total routing distance of the sprayers. The resulting mixed integer programming model that we develop is hard to solve using a commercial solver, owing to the dependencies caused by the spatio-temporal synchronization between the two operators. To solve large-scale instances effectively, we present an adaptive large neighborhood search metaheuristic that uses an intensive local search mechanism. We conduct an extensive computational analysis to assess the effectiveness of our solution approach and gain managerial insights into the problem by analyzing various models. The proposed metaheuristic yields high-quality solutions quickly, with an overall average improvement of 5.61\% over what is implemented in practice implying significant savings in time and cost.
}%

\noindent

\textit{\textbf{Keywords:}}
Vehicle routing problem with multiple synchronization constraints; vehicle routing; Vehicle routing problems with multiple synchronization constraints; synchronization; metaheuristic 


%


\section{Introduction}
\label{sec:intro}
In large farm spraying operations, two operators are needed, namely a fleet of sprayers and a tender tanker. The sole task of a sprayer is to apply fertilizers at infected areas in a farm while the sole task of the tender tanker is to refill a sprayer's tank whenever a sprayer needs a refill. Proper synchronization between sprayers and the tender tanker is crucial to achieving higher operational efficiency, better farming practices, and lower costs. While poor synchronization between the sprayers and the tender tanker will lead to higher waiting times for the sprayers to get a refill, an unnecessary increase in the number of refills for sprayers, an increase in the traveled distance of the tender tanker across the farm negatively affecting the soil compaction, and an increase to the operational cost in the form of fuel cost. Furthermore, the tender tanker itself has a limited tank capacity and it also performs a refilling by going back to the storage unit located at the farm's depot to get a refill. The sprayers and tender-tanker routing with synchronization and multiple trips problem hence considers the trade-off between increasing the frequency of refills and waiting time for sprayers to get a refill as well as the optimization of routing for the sprayers and the tender tanker. The goal of the sprayers-tender tanker routing synchronization problem with multiple trips is to determine the optimal assignment of infected locations of a farm to be served by a fleet of sprayers, the routing of sprayers, the refilling locations for the sprayers to get a refill, the scheduling of the refills, and the routing plan of the tender tanker with the possibility of performing multiple trips. \\

The problem we study in this paper shares characteristics and features with the well-known vehicle routing problem (VRP) with multiple synchronization constraints studied in the literature, see \cite{soares2023synchronisation} for a review. The VRP with multiple synchronization constraints can be used to model many real-world applications such as staff scheduling \cite{mankowska2014home}, home care delivery \cite{hashemi2020vehicle}, military aircraft mission planning \cite{quttineh2013military}, and two-echelon distribution systems \cite{grangier2016adaptive}. The practical importance of the VRP with multiple synchronization constraints arises from the fact that a wide range of variations in operating rules and constraints encountered in real-life applications can be modeled as VRP with multiple synchronization constraints.\\

In this paper, we investigate improving the operational efficiency of spraying in large farms by presenting an extension of the VRP with multiple synchronization constraints with multiple trips and realistic assumptions. In this problem, a single vehicle performs multiple trips to the depot and is responsible for serving a fleet of vehicles, on the other hand, the task of the fleet of vehicles is to apply fertilizer at infected locations in a farm to increase the agriculture yield and ensure healthier crops. In this context, synchronization is both spatial and temporal in the sense that the decision maker needs to decide on when and where each refilling operation takes place (recall that a refilling operation calls for synchronization between the tender tanker and a sprayer). Note that by deciding on the locations and scheduling of refillings, the complexity of the problem increases significantly.\\ 

There is little literature on VRPs with multiple synchronization constraints considering spatial and temporal aspects; to the best of our knowledge, this is the first study of this type of problem with multiple trips. The research we present in this paper has been inspired by a problem faced by a company providing spraying services to farmers. Farm spraying planning has to decide which sprayer is assigned to an infected area in a farm to be sprayed, the sequence of areas to be visited by each sprayer, the refilling scheduling, and the routing of the tender tanker. The tender tanker is a heavy tank and reducing its travel distances is very important to farmers who are concerned about soil integrity. Similarly, reducing the traveling time of the sprayers and the waiting time to get a refill allows them to perform more spraying, which in the long run increases their productivity.

Our work aims to develop and solve a mathematical model that will assist companies performing spraying services to achieve greater efficiency by reducing the routing and waiting times of sprayers and reducing the routing time of the tender tanker to preserve soil integrity. We present a summary of our contributions as follows:
\begin{itemize}
    \item[-] To the best of our knowledge and according to our reviewed literature in Section \ref{sec:litRev}, our paper is the first to propose and analyze a problem similar to the sprayer-tanker synchronized assignment routing problem with multiple trips feature.
    \item[-] Because of the challenges of solving large-scale instances of the mathematical model using a commercial solver, we propose an adaptive large neighborhood search (ALNS) metaheuristic to address the problem and provide high-quality solutions to large-scale practical size instances. The metaheuristic combines an effective heuristic, which yields an initial feasible solution in a short computational time and solves a mathematical model.
    \item[-] We build comprehensive computational sets of experiments using randomly generated instances to evaluate the efficacy of our approach and gain insights. Our results show that the suggested metaheuristic can solve large instances in a short period of time.
    \item[-] Additionally, our meta-heuristic can solve large-scale instances within a fraction of the time used by the solver.
    \item[-] We demonstrate that our approach yields an overall improvement of 5.61\% compared to the approach implemented in practice.
\end{itemize}

The rest of this paper is organized as follows. Section \ref{sec:litRev} introduces the relevant literature. Section \ref{sec:Operation} presents a clear of the problem we study. A mathematical programming formulation for the problem is developed in Section \ref{sec:mathModel}. Section \ref{sec:ALNS} proposes an ALNS metaheuristic to solve the problem. Computational results are reported in Section \ref{sec:Computational}. Finally, conclusions are presented in Section \ref{sec:conclude}.

\section{Literature Review}
\label{sec:litRev}
The problem we study in this work is a combination of the Vehicle Routing Problem with Multiple Synchronization constraints (VRPMSs) and the VRP with multiple trips.  According to \citet{drexl2012synchronization}, VRPMSs can be defined as ``A VRPMS is a vehicle routing problem in which more than one vehicle may or must be used to fulfill a task." What sets the VRPMSs apart from their VRPs counterparts is the nature of synchronization provided. As such, \cite{soares2023synchronisation} provides four categories of VRPMSs. These categories are routing with synchronization of schedules (e.g., \cite{fedtke2017gantry}), routing with transfers or cross-docking requirements (e.g., \cite{sacramento2019adaptive, schermer2019matheuristic}), routing with trailers or passive vehicles (e.g., \cite{meisel2014synchronized}), and routing with autonomous vehicles (e.g., \cite{murray2015flying}). Our problem fits in the last category. Therefore, in this section, we provide a brief review of the state-of-the-art studies in this category.\\

Motivated by last-mile delivery in logistics operations, \cite{murray2015flying} is the first to introduce optimal routing and scheduling of a single drone and a single delivery truck. In this delivery system, a customer can be served by a truck or a drone. The authors present mixed integer linear programming formulations and heuristic solution approaches to solve problems of practical size. The main advantage of using drones in delivery is to reduce routing time and add flexibility to the system. In their model, however, a drone can only visit one location before getting back to the truck for a recharge or get a new package. Since then, there has been a large number of studies in the literature optimizing the routing and synchronization of a fleet of drones and a fleet of vehicles. For example, the study of \cite{schermer2019hybrid} considers m-Truck and n-Drone system, \cite{thomas2023collaborative} considers 1-Truck and n-Drone system, and \cite{najy2023collaborative} consider 1-Truck and 1-Drone system within inventory routing setup with a time-horizon. It is worth mentioning that the number of studies considering multiple trucks and multiple drone systems is substantially less than the studies considering one truck and one drone or the one truck multiple drones systems. Interested readers are referred to \cite{boysen2021last} for a literature review on last-mile delivery with drones or robots assisting trucks. \\

\cite{boysen2018scheduling} propose an innovative technique for last-mile delivery using robots and a truck. The problem is concerned with scheduling procedures that determine the truck route along robot depots and drop-off points where robots are launched, such that the weighted number of late customer deliveries is minimized. In their work, the drop-off locations of the drones are a predefined set of locations. \\

\cite{sacramento2019adaptive} study multiple-truck and drone delivery system with the objective of minimizing delivery costs. Due to the difficulty of solving large instances to optimality, the authors propose an ALNS metaheuristic. The authors present the tests to demonstrate the benefits of including a drone within the delivery system. \\

\cite{li2020two} introduce the two-echelon vehicle routing problem with time windows and mobile satellites that optimizes delivery routes for a fleet of truck-drone combinations. In the first echelon, parcels are delivered using vans from a distribution center to customers and in the second echelon drones are dispatched from mobile-satellite vans to serve customers. The authors propose an adaptive large neighborhood search heuristic as an algorithm to solve the developed mathematical model. \\

\cite{coindreau2021parcel} propose a mixed-integer linear programming formulation and an adaptive large neighborhood search (ALNS) for the drone-truck routing problem motivated by the case of a large European logistics provider. Their model considers realistic delivery features such as time windows, limited drone autonomy, and the eligibility of clients to be served by drones. The objective function in their model aims at minimizing the fixed daily vehicle fares, driver wages, and fuel and electricity consumption to power trucks and drones.\\

\cite{tamke2021branch} study multi-truck and multi-drone delivery for delivery with the goal of minimizing the make-span. They consider a homogeneous fleet of vehicles and each vehicle is equipped with one drone. Each drone can deliver one package before getting back to the truck for recharging/reloading. The authors introduced new valid inequalities and a separation routine for extended subtour elimination constraints and they proposed a branch-and-cut algorithm to solve large-scale instances.\\

\cite{ostermeier2023multi} formulate a multi-vehicle truck-and-robot routing problem where a truck might be loaded with more than one robot for delivery systems with the objective of minimizing the total routing cost. The authors developed a tailored heuristic solution approach based on a novel neighborhood search. In their system, a robot can serve at most one customer before returning back to the truck.\\

\cite{momeni2023new} studied a vehicle routing problem with drones with the objective of minimizing total delivery time and with the ability to deliver postal packages in the path of drones at different heights. Their model allows the drone to deliver more than two packages per trip and the drone's energy consumption is computed by taking into account wind speed, the weight of the postal parcel, and the weight of the drone's body. Their model assumes one truck and more than one drone. The authors proposed a two-phase algorithm based on the nearest neighborhood and local search is presented to solve the developed mathematical model.\\

As far as studies in the literature considering VRP with multi-trips, \cite{mingozzi2013exact} describe two set partitioning formulations of the multiple-trip VRP and study valid lower bounds, based on the linear relaxations of both formulations enforced with valid inequalities, that are embedded into an exact solution method. \cite{cattaruzza2016multi} study multiple-trip vehicle routing problem with time windows and release dates and propose an efficient labeling algorithm to solve the new problem. \cite{franccois2016large} consider the multi-trip vehicle routing problem arising when customers are close to each other or when their demands are large. They design two large neighborhood search heuristics and compare them against solving the problem by combining vehicle routing heuristics with bin packing routines in order to assign routes to vehicles. In their work, \cite{franccois2016large} also provide insights into the configuration of the proposed algorithms by analyzing the behavior of several of their components. \cite{franccois2019adaptive} study multi-trip VRP with time windows and developed an ALNS as a solution approach. Unlike previous heuristic solution methods for multi-trip vehicle routing problems which separate routing and assignment phases, their work presented a new integrated approach that outperforms the previous methods. \cite{cattaruzza2018vehicle} present an exhaustive survey and unified view on multi-trip VRP mathematical formulations and the solution approaches. \\



Despite the aforementioned similarity between the cited work and our work, there are also key differences: (1) all studies considering vehicle routing setup with truck and drone systems assume that a node can be served by a truck or a drone, on the other hand, in our setup, the roles of the sprayers and the tender tanker can not be exchanged, (2) for truck-drone delivery system, each drone is assigned to a truck, by contrast, in our work, the tender tanker serves all sprayers in the field, and (3) in all studies considering vehicle routing problem with multiple synchronization constraints, no study incorporated the multi-trip feature. This paper is the first to address the vehicle routing problem with multiple synchronization constraints and multiple trips considering the features of a fleet of sprayers and a tender tanker spraying operations. Despite the practical importance of this problem for agriculture companies and farmers, there is an evident gap in the literature that inspired us to develop a new model considering the features of the sprayers and tender tanker with synchronized refilling operations and multiple trips.

\section{Operation description}
\label{sec:Operation}
It is crucial in the farming business to efficiently manage a fleet of sprayers who perform spraying activities for the planted crops. We model a farm, for the purposes of spraying activities, as a graph with a set of locations/nodes and arcs connecting these nodes. An infected area is the coordinate where harmful plants exist (such as weeds) and it needs to be sprayed so that it does not grow and cause harm to the crop. In this work, we refer to an area infected with harmful plants to be infected as a \textit{``node"}. To properly limit agricultural damage caused by harmful plants, each infected location must be treated with a certain quantity of fertilizer, and these spraying operations are carried out by a fleet of sprayers. Typically, the type of sprayer used is the boom sprayer. This sprayer may be mounted on a tractor. Nozzles may be positioned and adjusted to provide adequate spray coverage. The boom sprayers must be driven through every row or every other row to get enough coverage.\\

The sprayer's tank is filled with fertilizer to perform the spraying operations. A typical size of a boom sprayer tank is 1,000 - 1,500 gallons. The sprayers move between the nodes and perform the spraying operations as instructed. Once a sprayer runs out of fertilizer, the tender tanker arrives at its node and refills its tank. On the other hand, the tender tanker's main role is to supply the sprayers with fertilizer by refilling their tank whenever they need a refill of fertilizer. A tender tanker's capacity is usually 10,000 - 15,000 gallons. \\

In case a sprayer runs out of fertilizer and a refill is needed while the tanker has not arrived at its location, the sprayer will wait for the tanker to arrive. Such waiting time is not desired as it is considered to be a loss in productivity and a waste of the sprayer's time. The decision of which sprayer should be assigned to which node, as well as the scheduling decisions of when and where the tanker should refill each sprayer whenever they run out of fertilizer throughout the operation, is a difficult task. This must be coupled with the routing of each sprayer and the tender tanker each day to serve the set of nodes to which they have been assigned. This is a complex integrated synchronized vehicle routing problem in which routing, scheduling, and assignment decisions should be addressed simultaneously, as well as dependencies between the tanker and the sprayer.\\

Inspired by a real-world problem faced by a company providing spraying services to farmers, we study and model the Vehicle Routing Problem with Multiple Synchronization Constraints and Multiple Trips (VRPMSC-MT) that we described above. The efficiency of spraying operations for our partner's company is based on minimizing the spraying costs that are based on the total time the sprayers spend from the moment they leave the depot till they get back to the depot.\\

\begin{figure}[t!]
\centering
\includegraphics[width=12cm]{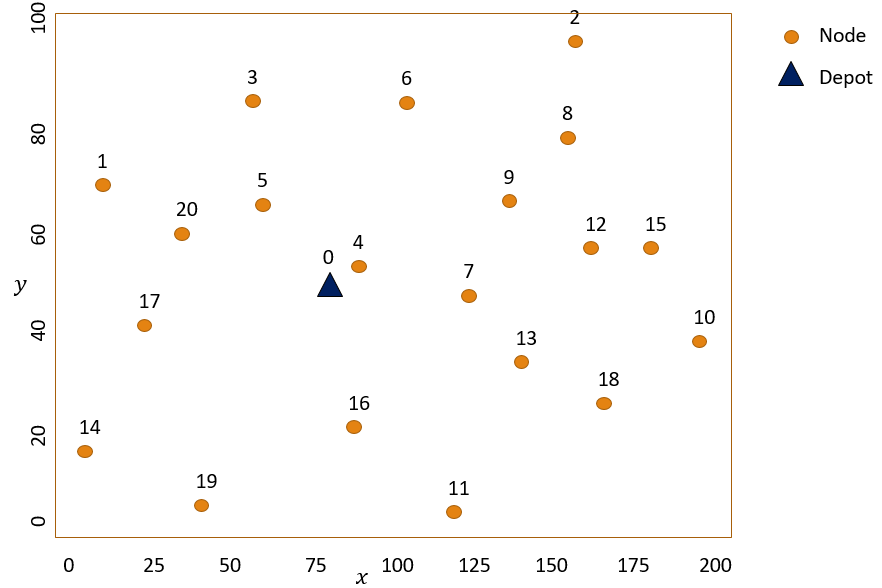}
\caption{Set of nodes/locations to be sprayed}
\label{fig:Nodes}
\end{figure}

We now provide a small example of the operation along with the added value of the optimization approach. In this small example, we assume that there are 20 nodes/locations to be sprayed on a farm. Nodes 1-20 are the set of nodes that need to be sprayed and node 0 represents the depot. Figure \ref{fig:Nodes} illustrates a set of nodes that need to be sprayed. The current practice our partner company implements is based on a two-stage optimization framework known as the route-first cluster-second. In the first stage, they run a capacitated vehicle routing algorithm to find the shortest routes connecting all nodes. In that capacitated vehicle routing model they set the capacity of a vehicle (i.e., sprayer) to the same capacity as the sprayer's tank. For illustration purposes, the sprayer capacity was set to 8 gallons, the tender tanker capacity to 30 gallons, and the demand for each node was set to 2 gallons.\\

\begin{figure}[t!]
\centering
\includegraphics[width=12cm]{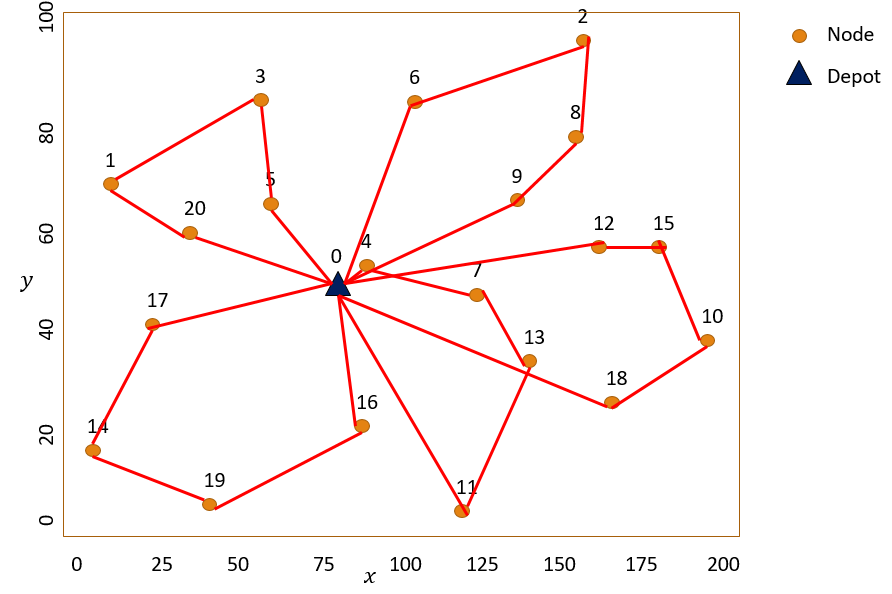}
\caption{Routes of capacitated vehicle routing problem}
\label{fig:VRP..}
\end{figure}

\begin{figure}[t!]
\centering
\includegraphics[width = 14cm,  height = 7cm]{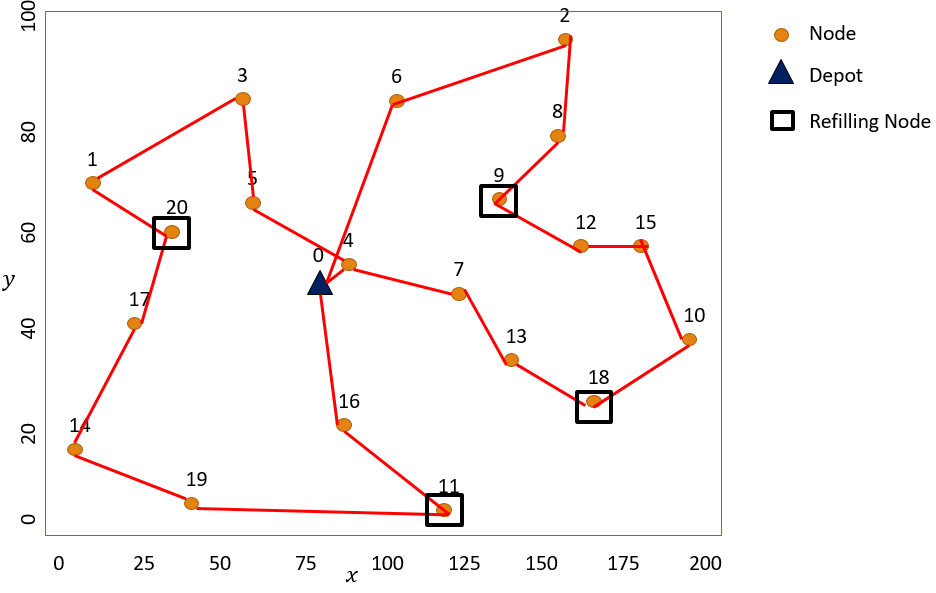}
\caption{Routes of the sprayers and the tender tanker: two-stage framework}
\label{fig:ST1}
\end{figure}

Applying a capacitated VRP model to find the shortest routes serving all the nodes while respecting the sprayer's tank capacity results in generating five routes. Each route starts and ends at the depot see Figure \ref{fig:VRP..}. The generation of these routes completes stage 1. The next steps are to identify the refilling nodes and the routes of the tender tanker. Identifying the refilling nodes is fairly simple now, it is based on the fact that the end of each route represents a point where the sprayer runs out of fertilizers; hence, the last node of each route is the refilling point. Once the sprayer's tank is filled again, routes can be connected and spraying can continue. The refilling nodes are then connected to form big routes. The number of generated big routes equals the number of sprayers available to spray the set of nodes. Note that the number of sprayers dedicated to spray is determined by the farm manager. For instance in Figure \ref{fig:VRP..} we identify node 20 as a refill node for the sprayer since it is the end of a route. We then connect this node to node 17 in the next route, see Figure \ref{fig:ST1}A.\\

If a sprayer requires a refill after completing spraying operations while the tender tanker is not there, the sprayer will wait for the tender tanker to arrive to start the refilling process. Such a policy yields a feasible solution to the problem. Applying this logic, the tender tanker's route can now be determined with the nodes of refills established. It is worth mentioning that the company we partner with prefers a consistent and clear policy regarding the timing of a refill to be communicated with the sprayers, namely either a refill might take place at a node before the sprayer starts spraying or after spraying is completed at a given node. Therefore, we assume that refilling happens only after a sprayer completes spraying a node/location.\\

\begin{figure}[t!]
\centering
\includegraphics[width = 12cm,]{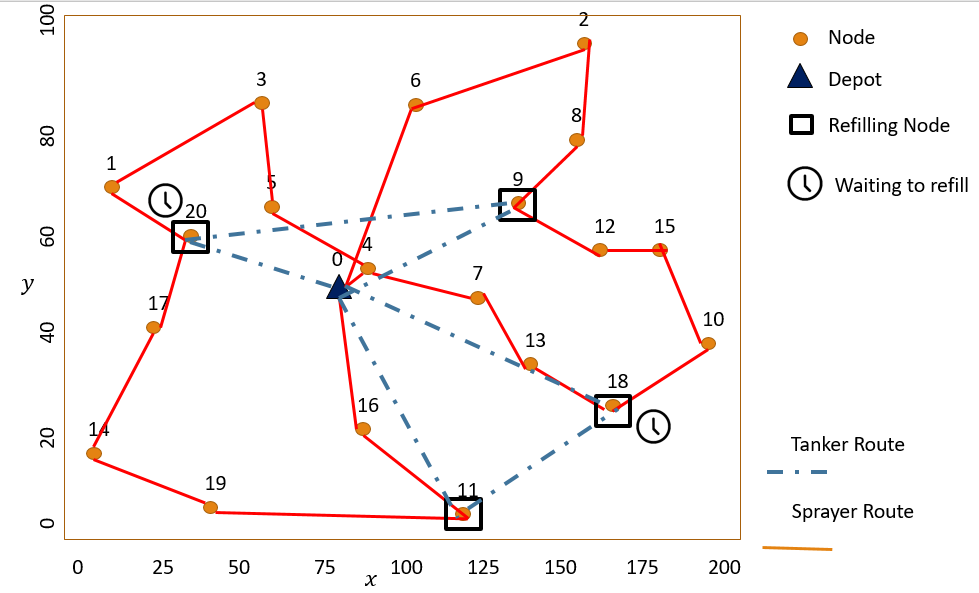}
\caption{Feasible solution from the two-stage optimization framework}
\label{fig:illus}
\end{figure}

The tender tanker route assignment is done manually, for instance as \textit{sprayer1} finishes spraying at node 9, the tender tanker will refill \textit{sprayer1}. \textit{Sprayer1} continues to the next node 12 and the tender tanker also proceeds to node 20 to refill \textit{sprayer2}. Figure \ref{fig:ST1}B shows the tanker's route. As a result, the tender tanker's route decisions are now established. See Figure \ref{fig:illus} for the full illustration of the example. \\

The performance metrics are the total routing time and the waiting time of the sprayers. The overall routing time and waiting time of the sprayers are calculated to be 72.38 and 1.60 units of time, respectively. Decomposing this integrated problem into two stages may not yield the optimal solution.\\

We solve the same instance using the optimization approach we present in the next section for the minimization of the total routing time and waiting time of the sprayers, see Section \ref{sec:model}. The total routing time under the optimization approach is 62.47 units of time and waiting time is zero. Clearly, the optimization approach we propose in this study is more effective in terms of decreasing the overall time spent by the sprayers when compared to the two-stage optimization framework. Hence, solving the problem using the optimization approach yields significant savings for farmers without compromising any performance measure.\\

\section{Synchronized Spraying Routing Problem}
\label{sec:mathModel}
\subsection{Mathematical model}
\label{sec:model}
The spraying planning involving a fleet of sprayers and a tender tanker problem is defined on a complete, undirected graph $G = (\mathcal{N}, \mathcal{A})$, where set $\mathcal{N}=\lbrace 0\cup \mathcal{N}_f\rbrace$ is comprised of the depot $D = \lbrace 0 \rbrace$ and nodes $\mathcal{N}_f=\lbrace 1,2,\cdots, N \rbrace$ to be sprayed and $\mathcal{A}$ is the set of undirected arcs. Denote the set of available trips as $\mathcal{K}=\lbrace 1,...,K\rbrace$ with members $k$ and $r$. Lastly, define a dummy copy of the depot as $N+1$,\\

The size of the fleet of sprayers is denoted as $NumSp$. Each node to be sprayed $i\in\mathcal{N}_f$ is associated with a certain quantity of fertilizer $q_i$ to be applied and a spraying (service) time defined as $s_i$. The spraying time parameter is calculated simply based on the spraying rate of the sprayer and the required spraying quantity. The working hours, which are considered as the planning horizon, are denoted by $tMax$. We assume that the network is fully connected and hence each pair of nodes in the network are connected via a direct arc. To each arc $(i, j) \in \mathcal{A}$ a travel time $t_{ij}$ is attributed. There is a fleet of homogeneous sprayers and a single tender tanker that are fully loaded and located at the central depot. We define a trip $k$ in a set of trips $\mathcal{K}$ for the tender tanker\footnote{A trip is defined as a route of the tender tanker starting at the depot to serve sprayers and ending at the depot.}. The capacity of the sprayer tank is $Q_s$ and the tender tanker's tank has a capacity of $Q_t$ and $Q_t > Q_s$. Finally, we define $M$ as a large number. The minimum number of trips that the tender tanker needs to perform can be calculated as $\left \lceil \frac{\sum_{i\in\mathcal{N}_f}q_i}{Q_t} \right \rceil$; however, setting the number of trips to a low value may result in less flexibility. As such, we have $\vert \mathcal{K}\vert = 2*\left \lceil \frac{\sum_{i\in\mathcal{N}_f}q_i}{Q_t} \right \rceil$ to balance model complexity due to the number of binary variables and the flexibility in modeling. \\

We assume that the refilling time needed to refill a sprayer is $\xi$ units of time regardless of the level of fertilizer in the sprayer's tank. This is a mild assumption since refilling involves three steps: 1. Setting up the hose to be attached to the sprayer's tank, 2. Supplying the sprayer's tank with fertilizer itself, and 3. De-attaching the hose from the sprayer's tank and getting back the equipment to the tender tanker. As such, supplying the sprayer's tank with fertilizer process is only one part of the refilling process. Likewise, the refilling time of the tender tanker at the depot is denoted as $\gamma$ and it is fixed regardless of the quantity of fertilizer that the tender tanker fills in its tank. \\

\subsubsection{Mathematical model formulation \label{mathmodel}}
To model the problem, the following decision variables are defined:

\begin{align*}
  x_{ij} &=
  \begin{aligned}[t]
    &\begin{cases}
      1, & \text{if arc $(i, j) \in \mathcal{A}$ is traversed by a sprayer,}\\
      0, & \text{otherwise,}
    \end{cases}
  \end{aligned} \\
  g_{ijk} &=
  \begin{aligned}[t]
    &\begin{cases}
      1, & \text{if arc $(i, j) \in \mathcal{A}$ is traversed by the tender tanker on trip $k\in\mathcal{K}$,}\\
      0, & \text{otherwise,}
    \end{cases}
  \end{aligned} \\
  \delta_{i} &=
  \begin{aligned}[t]
    &\begin{cases}
      1, & \text{if refilling a sprayer takes place at node $i \in \mathcal{N}_f$,}\\
      0, & \text{otherwise,}
    \end{cases}
  \end{aligned} \\
  z_{kr} &=
  \begin{aligned}[t]
    &\begin{cases}
      1, & \text{if trip $k \in \mathcal{K}$ proceeds trip $r \in \mathcal{K}$,}\\
      0, & \text{otherwise,}
    \end{cases}
  \end{aligned}\\
  \theta_{ik} &
  \begin{aligned}[t]
      & \text{  the arrival time of the tender tanker at node $i \in \mathcal{N}_f$ on the $k^{th}$ trip,}
  \end{aligned}\\
  y_{i} &
  \begin{aligned}[t]
       & \text{  the time of arrival of a sprayer at node $i \in \mathcal{N}_f$,}
  \end{aligned}\\
  a_{i} &
  \begin{aligned}[t]
       & \text{  the time of spraying service start time by a sprayer at node $i \in \mathcal{N}_f$,}
  \end{aligned}\\
  w_{ik} &
  \begin{aligned}[t]
       & \text{  the time of refilling start at node $i \in \mathcal{N}_f$ on the tender tanker $k^{th}$ trip,}
  \end{aligned}\\
  h_{ik} &
  \begin{aligned}[t]
       & \text{ remaining quantity of fertilizer in the tender tanker's tank upon arrival to node $i \in \mathcal{N}$ on the $k^{th}$ trip,}
  \end{aligned}\\
  l_i & 
  \begin{aligned}[t]
      & \text{  remaining quantity of fertilizer in the sprayer's tank upon arrival to node $i \in \mathcal{N}_f$,}
  \end{aligned}\\
    m_i & 
  \begin{aligned}[t]
       & \text{ sprayer's waiting time at node $i \in \mathcal{N}_f$.}
  \end{aligned}\\
  v_i & 
  \begin{aligned}[t]
       & \text{  quantity of refill to a sprayer at node $i \in \mathcal{N}_f$.}
  \end{aligned}
\end{align*}

Table~\ref{Table:Notation} summarizes all notation. Using this notation, we formulate the following MILP:
\begin{table}[h!]
\small
  \renewcommand{\arraystretch}{0.4}
\caption{Notation}
\label{Table:Notation}
\begin{tabular}{ll}
\hline 
\multicolumn{2}{l}{\textbf{Sets}}                       \\ \hline
$\mathcal{N}$         & Set of nodes                        \\
$\mathcal{N}_f$         & Set of nodes to be sprayed                        \\
$\mathcal{K}$         & Set of tanker's trips                        \\
\hline
\multicolumn{2}{l}{\textbf{Indices}}                       \\ \hline
$i,j$         & Node indices                        \\
$k,r$           & Trip indeces                        \\ \hline
\multicolumn{2}{l}{\textbf{Parameters}}                        \\ \hline
$numSp$           & Number of sprayers  \\
$q_i$             & Quantity of fertilizer to be applied at node $i\in\mathcal{N}_f$          \\
$s_i$             & Service time of node $i\in\mathcal{N}_f$          \\
$tMax$            & Working hours/length of the planning horizon          \\
$t_{ij}$    & Travelling time between nodes $i$ and $j$      \\
$Q_s$      & Capacity of the sprayer \\
$Q_t$      & Capacity of the tanker \\
$\xi$      & Time needed to refill a sprayer \\
$\gamma$      & Time needed to refill the tanker at the depot \\
$M$       & Big positive number \\ \hline
\multicolumn{2}{l}{\textbf{Decision variables}} \\ \hline
$x_{ij} $ & Binary variable that is equal to 1 if and only if sprayer travels along arc $(i,j)$ \\
$g_{ij}^k$     & Binary variable that is equal to 1 if and only if the tanker travels arc $(i,j)$ in trip $k$ \\
$\delta_i$     & Binary variable that is equal to 1 if and only if there is a refill at node $i\in\mathcal{N}_f$ \\
$z_{kr}$     & Binary variable that is equal to 1 if and only if trip $k\in\mathcal{K}$ proceeds trip $r\in\mathcal{K}$ \\
$\theta_i^k$     & Arrival time of the tender tanker at node $i\in\mathcal{N}_f$ on the $k^{th}$ trip\\
$y_i$     & Time of arrival of a sprayer at node $i\in\mathcal{N}_f$ \\
$a_i$     & Time of spraying service start time by a sprayer at node $i\in\mathcal{N}_f$ \\
$w_i^k$     & Refilling start time of the tender tanker at node $i\in\mathcal{N}_f$ on the $k^{th}$ trip\\
$h_i^k$     & Level of fertilizer of the tender tanker at arrival to node $i\in\mathcal{N}_f$ on the $k^{th}$ trip\\
$l_i$       & Level of fertilizer of sprayer when arriving to node $i\in\mathcal{N}_f$ \\
$m_i$       & Waiting time of a sprayer at node $i\in\mathcal{N}_f$ \\
$v_i$       & Quantity of fertilizer supplied to a sprayer at node $i\in\mathcal{N}_f$ \\
 \hline
\end{tabular}
\end{table}
 
\newpage 
\begin{eqnarray}
\label{eq:FObjt}
\displaystyle{\min}&& \displaystyle{\sum_{i\in \mathcal{N}}\sum_{j\in \mathcal{N}: i\neq j} t_{ij}x_{ij}} + \sum_{i\in \mathcal{N}_f} m_{i}  + \sum_{i\in \mathcal{N}_f} \xi \delta_{i}\\
s.t.
\label{eq:FOne1}
& &\displaystyle{\sum_{i\in \mathcal{N}:i\neq j} x_{ij} = 1 \quad \forall \quad j\in \mathcal{N}_f,} \\
\label{eq:FThree2}
& & \displaystyle{\sum_{i\in\mathcal{N}_f:i\neq j} x_{ij} = \sum_{i\in\mathcal{N}_f:i\neq j} x_{ji} \quad \forall \quad } j\in \mathcal{N}, \\
\label{eq:FTNumSp}
& & \displaystyle{\sum_{i\in\mathcal{N}_f:i} x_{0i} =\sum_{i\in\mathcal{N}_f:i} x_{i0} = numSp,} \\
\label{eq:FFour3}
& & \displaystyle{\sum_{i\in\mathcal{N}_f:i\neq j} g_{ijk} = \sum_{i\in\mathcal{N}_f:i\neq j} g_{jik} \quad \forall \quad  j\in \mathcal{N}, k \in \mathcal{K},}\\
\label{eq:FFive44}
& & \displaystyle{\delta_{i} = \sum_{j\in \mathcal{N}:j\neq i}\sum_{k\in \mathcal{K}} g_{jik} \quad \forall \quad i\in\mathcal{N}_f,}  \\
\label{eq:FSix5}
& & \displaystyle{ \vert \mathcal{K}\vert-\sum_{k\in\mathcal{K}} \sum_{r\in\mathcal{K}:k<r} z_{kr} = 1,}\\
\label{eq:FFour25}
& & \displaystyle{\sum_{r\in\mathcal{K}:r<k} z_{rk} \leq 1 \quad \forall \quad s\in\mathcal{K},}\\
\label{eq:FFour26}
& & \displaystyle{\sum_{r\in\mathcal{K}:k<r} z_{kr} \leq 1 \quad \forall \quad k\in\mathcal{K},}\\
\label{eq:FEight7}
&&\displaystyle{y_i + s_i  \leq w_{ik} \quad \forall i\in \mathcal{N}_f, k \in \mathcal {K},}\\
\label{eq:FTen14}
& & \displaystyle{\theta_{ik}\leq w_{ik}  \quad \forall k \in \mathcal {K},}  \\
\label{eq:FTen9}
& & \displaystyle{y_j\geq y_i+(s_i + t_{ij})x_{ij} + \xi \delta_i + m_i - M(1-x_{ij}) \quad \forall \quad i,j\in \mathcal{N}_f: i\neq j,}\\
\label{eq:FTen91}
& & \displaystyle{y_j\leq y_i + (s_i + t_{ij})x_{ij} + \xi  \delta_i + m_i + M(1-x_{ij}) \quad \forall \quad i,j\in \mathcal{N}_f: i\neq j,}\\
\label{eq:FTen10}
& & \displaystyle{y_i\geq t_{0i}x_{0i} \quad \forall \quad i\in \mathcal{N}_f,}\\
\label{eq:FTen25}
& & \displaystyle{w_{ik}+t_{i0}g_{i0k}+\xi\delta_i \leq \theta_{N+1,k}+M(1-g_{i0k}) \quad \forall \quad i\in \mathcal{N}_f, k \in \mathcal {K},}  \\
\label{eq:FTen177}
& & \displaystyle{\theta_{N+1,k} + \gamma \leq \theta_{0,r} + M(1-z_{kr})  \quad \forall \quad  k,r \in \mathcal {K}:k<r,}\\
\label{eq:FTen12}
& & \displaystyle{s_i + y_i + t_{i0}x_{i0}\leq tmax  \quad \forall \quad i\in \mathcal{N}_f}\\
\label{eq:FTen13}
& & \displaystyle{\theta_{ik}+t_{i0}g_{i0k}\leq tmax  \quad \forall k\in\mathcal {K},}  \\
\label{eq:FTen155}
& & \displaystyle{w_{ik}+(t_{ij}+\xi)g_{ijk} - M(1-g_{ijk})\leq \theta_{jk} \quad \forall \quad i,j\in \mathcal{N}_f:i\neq j, \quad k \in \mathcal {K},}\\
\label{eq:FTen22}
& & \displaystyle{m_i\geq \theta_{ik} - (y_i +s_i) - M(1-\theta_{ik})\quad \forall \quad i\in \mathcal{N}_f \quad k\in \mathcal{K},}\\
\label{eq:FTen15}
& & \displaystyle{l_j \leq l_i - q_ix_{ij} + v_i +  M(1-x_{ij}) \quad \forall  j\in \mathcal{N}_f, i\in\mathcal{N}_f:i\neq j}\\
\label{eq:FTen16}
& & \displaystyle{ 0 \leq h_{jk} \leq h_{ik} - v_i + M(1-g_{ijk}) \quad \forall  j\in \mathcal{N}_f, i\in\mathcal{N},  k\in \mathcal{K},}\\
\label{eq:FTen166}
& & \displaystyle{v_{i} \leq Q_s + Q_s(1-\delta_i) - l_i + \sum_{j\in\mathcal{N}:j\neq i}q_ix_{ji} \quad \forall  i\in \mathcal{N}_f,}\\
\label{eq:FTen167}
& & \displaystyle{v_{i} \leq Q_s\delta_i \quad \forall  i\in \mathcal{N}_f,}\\
\label{eq:FTen17}
& & \displaystyle{l_0 = Q_s,}  \\
\label{eq:FTen18}
& & \displaystyle{h_0 = Q_t,}  \\
\label{eq:FTen19}
& & \displaystyle{l_j\leq Q_s  \quad \forall \quad i\in \mathcal{N}_f,}  \\
\label{eq:FTen70}
& & \displaystyle{h_{jk}\leq Q_t  \quad \forall \quad i\in \mathcal{N}_f, k\in \mathcal{K},}  \\
\label{eq:FTen27}
& & \displaystyle{l_i\geq Q_s x_{0i}  \quad \forall \quad i\in \mathcal{N}_f,}\\
\label{eq:FIntg}
& &\displaystyle{g_{ijk}, x_{ij},\delta_{i} \in \lbrace 0,1\rbrace}, m_i,v_i\theta_{ik},y_{i},l_{i},w_{ik}, h_{ik} \geq 0 \quad \forall i,j,k. \
\end{eqnarray}

The objective function (\ref{eq:FObjt}) minimizes the sprayers' travel time, refilling time, and waiting time. The objective function is consistent with the preferences of the farmers, namely, complete the spraying operations with the minimum time wasted for sprayers being idle or needing a refill and minimize the distance traveled by the sprayers and the tanker since their movement in the farm/field negatively affects soil compaction. Note that we use the traveling time as a function to reflect the traveled distance since they are correlated and to avoid having a multi-objective optimization model. This is a mild assumption. The constraints can be divided into four groups as follows.\\

\textbf{Sprayers and tender tanker routing:} Constraints (\ref{eq:FOne1}) ensure that each node is visited once by a sprayer. Constraints (\ref{eq:FThree2}) ensure that if a node is visited by a sprayer; then the sprayer leaves that node. Similarly, constraints (\ref{eq:FFour3}) ensure that if a node is visited by the tender tanker; then the tender tanker leaves that node. Constraints (\ref{eq:FFive44}) ensure that a refill process may happen only if the tanker visits that node. Constraint (\ref{eq:FSix5}) implies that there isn only one tender tanker. Constraints (\ref{eq:FFour25})-(\ref{eq:FFour26}) prevent variables $z_{kr}$ from being artificially set to one. They are needed to make sure that the left-hand side of Constraint (\ref{eq:FSix5}) indeed represents the number of tender tankers used.

\textbf{Arrival and refilling scheduling:} Constraints (\ref{eq:FEight7}) and (\ref{eq:FTen14}) ensure that refilling starting time at node $i\in\mathcal{N}_f$ happens after the spraying is complete or the arrival of the tender tanker, whichever is later. Constraints (\ref{eq:FTen9})-(\ref{eq:FTen91}) determine the arrival time of the sprayer at each node $i\in\mathcal{N}_f$. Constraints (\ref{eq:FTen10})-(\ref{eq:FTen25}) determine the arrival time of sprayer and the tender tanker at each node $i\in\mathcal{N}_f$ immediately after leaving the depot. Constraints (\ref{eq:FTen177}) calculate the start and end of each trip by the tender tanker. Constraints (\ref{eq:FTen12}) and (\ref{eq:FTen13}) ensure that the maximum length of a trip does not violate the workload of a sprayer or the tender tanker. Constraints (\ref{eq:FTen155}) determine the arrival time of the tender tanker at each refilling node. Constraints (\ref{eq:FTen22}) calculate the wait time of the sprayer at node $i$.

\textbf{Quantity of fertilizer calculation:} Constraints (\ref{eq:FTen15})-(\ref{eq:FTen16}) determine the fertilizer level of both the sprayers and the tender tanker after visiting a node. Constraints (\ref{eq:FTen166}) and (\ref{eq:FTen167}) calculate the proper value of the refilling quantity for a sprayer at node $i$ based on the level of fertilizer in its tank. Constraints (\ref{eq:FTen17})-(\ref{eq:FTen18}) ensure that both the sprayer and the tender tanker leave the depot at full tank capacity. Constraints (\ref{eq:FTen19})-(\ref{eq:FTen70}) ensure that the level of fertilizer in the tanker and sprayers' tank at node $i$ does not exceed their tank capacity. Constraint (\ref{eq:FTen27}) ensures that the sprayer's level of fertilizer is at full capacity when it visits the first node $i$. Lastly, constraint (\ref{eq:FIntg}) imposes the domain of each decision variable.\\

Clearly, mathematical model (\ref{eq:FObjt})-(\ref{eq:FIntg}) aims at minimizing the traveling time, waiting time, and refilling times. We call that mathematical model $Model1$. We generate models to account for various perspectives by modifying $Model1$. The following model modifications are presented:

\begin{eqnarray}
\label{eq:ObjF}
\displaystyle{\min} && \displaystyle{\alpha}\\
\label{eq:FNew1}
&&\displaystyle{\alpha \geq  y_i + s_i + t_{i0}x_{i0} \quad \forall i\in \mathcal{N}_f,}
\end{eqnarray}

We first focus on minimizing the latest arrival of the sprayers. We replace objective function (\ref{eq:FObjt}) of $Model1$ by equation (\ref{eq:ObjF}) and add constraints (\ref{eq:FNew1}). We call that model $Model2$. \\

Lastly, $Model1$ is modified by minimizing the total routing time only and ignoring the refill and waiting times. We call that model $Model3$.

\section{Metaheuristic for Vehicle Routing Problem with Multiple Synchronization Constraints and Multiple Trips}
\label{sec:ALNS}
The problem we discuss in this study is an NP-hard problem since it is a VRP with multiple synchronization and multiple trips. The mathematical formulation we present in this paper can be solved using broad-sense solvers such as Gurobi only when the size of the instances is small. Given examples of large-scale real-world settings, the problem may be intractable. As a result, we develop a metaheuristic to find high-quality solutions for large instances in short time periods.\\

The metaheuristic we develop consists of two phases, the first phase aims at finding an initial feasible solution that respects all constraints and the second phase aims at improving the initial feasible solution iteratively using ALNS. 

\subsection{Initial feasible solution procedure}
\label{sec:feasible solution}
Ideally, we would like to develop a mathematical model that shares some similarities with the original mathematical model (\ref{eq:FObjt})-(\ref{eq:FIntg}) with some relaxations to get an initial feasible solution, see \cite{archetti2017matheuristic,alkaabneh2022multi}. However, given the complexities of model (\ref{eq:FObjt})-(\ref{eq:FIntg}) and the inter-dependencies between the tender tanker and the sprayers, such an approach is not practical from a computational time point of view. Instead, we need a fast method to generate initial feasible solutions that will be later on improved using the ALNS and a mathematical model. That said, we opt to utilize a clustering algorithm to partition the nodes into clusters based on the $x-y$ coordinate of each node with a number of clusters equals the number of available sprayers and adapt a cluster-first-route-second procedure to establish the routes of the sprayers. The clustering algorithm we use in our work is based on \textit{k-}means unsupervised machine learning algorithm. Once each cluster is established, the routing of a sprayer assigned to that cluster is calculated aiming at minimizing the total traveled distance via a simple TSP model, then the refill nodes are established. Assuming that there is no waiting time for refilling process (i.e., assuming that a tender tanker is readily available to refill a sprayer) we calculate \textit{estimated arrival and departure times} at each node. We describe these arrival and departure times as \textit{estimated} since they do not take into account the tanker synchronization. Using the \textit{estimated departure times} of each of the nodes that will be refilled, the sequence of tender tanker is established by assuming that the tender tanker will visit the refill positions following an ascending order based on the \textit{estimated departure} times calculated in the previous step. Finally, we calculate the actual arrival and departure times of the sprayers and the tender tanker since we have the routes of each one of them. \\

However, such a solution may not be feasible to model (\ref{eq:FObjt})-(\ref{eq:FIntg}) due to violating constraints (\ref{eq:FTen12})-(\ref{eq:FTen13}) (i.e., a sprayer or the tanker does not finish its trip and get back to the depot within $tMax$). To ensure feasibility, we use a simple fix to reduce the length of the busiest route by assigning nodes from the busiest route to other routes. We elaborate more on that in the following paragraph using an example. \\

\begin{figure}[H]
\centering
\begin{subfigure}{1.0 \textwidth}
\centering
    \includegraphics[width=0.6\textwidth]{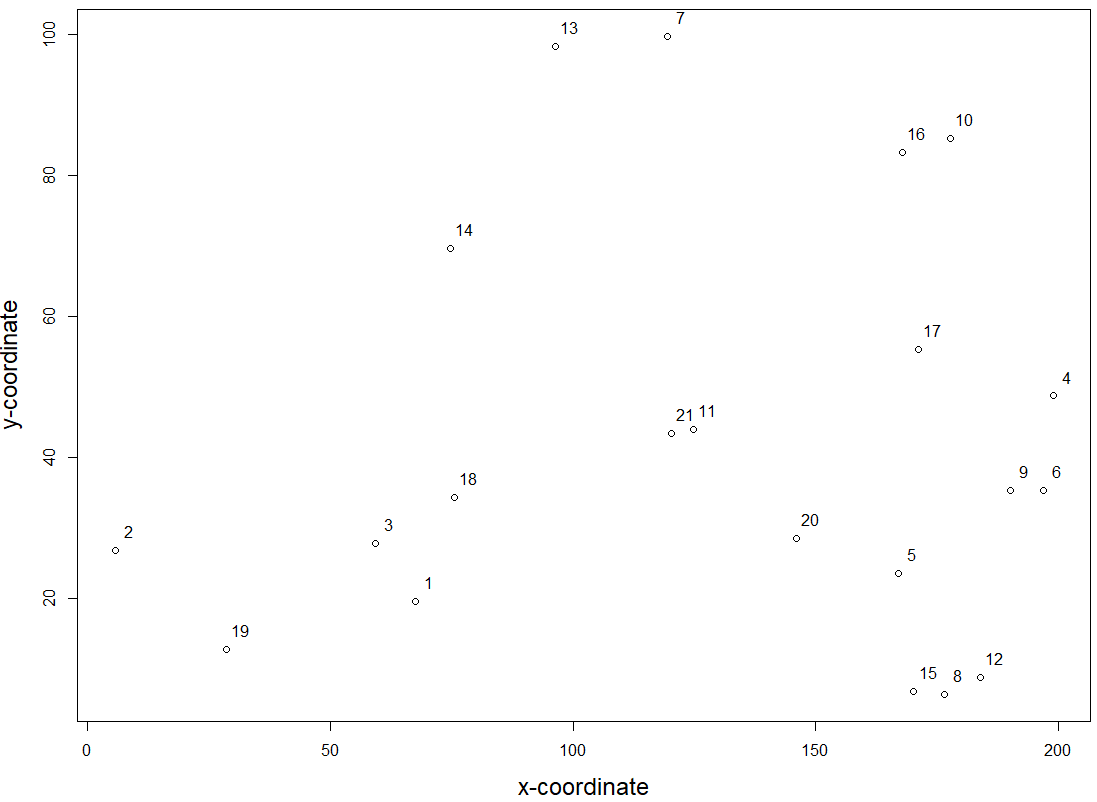}
    \caption{Map of nodes.}
    \label{fig:first}
\end{subfigure}

\hfill

\begin{subfigure}{1.0 \textwidth}
\centering
    \includegraphics[width=0.6\textwidth]{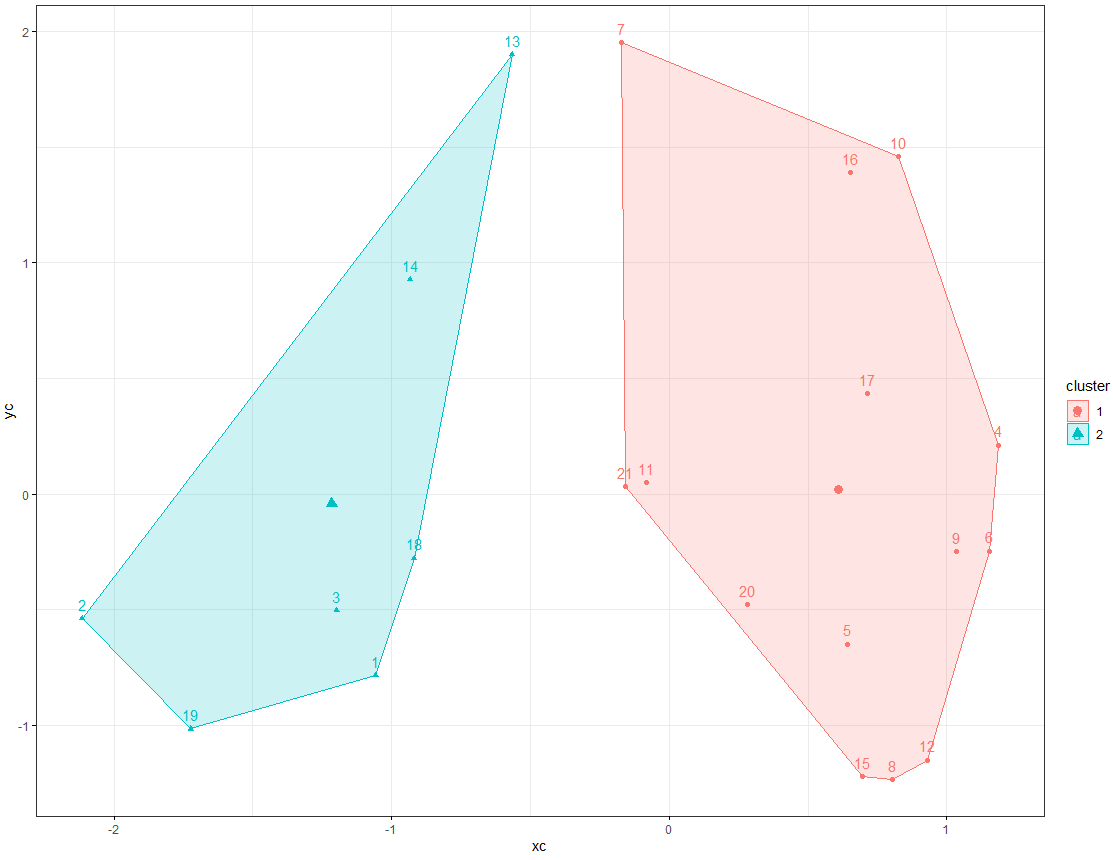}
    \caption{Clustering.}
    \label{fig:second}
\end{subfigure}
\caption{K-means clustering.}
\label{fig:feasibleSol}
\end{figure}

One of the downsides of using k-means clustering algorithms is that the resulting clusters may not be balanced in terms of the size of each cluster. For instance, given a set of 21 nodes displayed in Figure \ref{fig:first} and assuming the availability of two sprayers. Upon applying k-means clustering, node 1, 2, 3, 13, 14, and 19 belong to the same cluster, call it cluster 1, and the rest of the node belong to the cluster 2 as displayed in Figure \ref{fig:second}. Note that the route visiting the node in cluster 2 exceeds $tMax$ and hence the solution is not feasible. Our fixing procedure starts by searching for the points in cluster 2 that are close to other clusters, in this example cluster 1. We insert node 7 and 21 to cluster 1 and remove them from cluster 2. Then we complete establishing the feasible solution. If at the end of the procedure one of the clusters causes infeasible solution with respect to $tMax$, the fix process continues until a feasible solution is found. Once a feasible solution is found, the metaheuristic proceeds to phase 2 which is the ALNS. 

\subsection{Calculating refilling locations, waiting times, and trips of the tender tanker}
\label{sec:soldRep}
In our work, we start by constructing the routing variables of sprayers, then we calculate the refilling points and the waiting times of sprayers, and finally we establish the tender tanker's routing. As such, in this section, we provide the details of calculating these values. But first, let us introduce the representation of a solution for the VRPMSC-MT. \\

The routing of sprayers is represented by an ordered sequence:

\begin{align*}\label{eq:pareto mle2}
Route_{sp} &=  [\langle (str_0,end_0),l_0,y_0,m_0,v_0,waiting_{0} \rangle, \langle (str_1,end_1),l_1,y_1,m_1,v_1,waiting_{1} \rangle,...,\\
&\langle (str_{N-1},end_{N-1}),l_{N-1},y_{N-1},m_{N-1},v_{N-1},waiting_{N-1} \rangle,\langle (str_{N},end_{N}),l_{N},y_{N},m_{N},v_{N},waiting_{N} \rangle]
\end{align*}

for each tuple, $(str_i,end_i)$ represents an active arc indicating the traveling of a sprayer from node $str_i$ to node $end_i$, $l_i$ represents the level of fertilizer upon arrival to node $i$, $y_{i}$ represents the arrival time at node $i$, $m_i$ represents the refilling time at node $0$ if and only if refilling happens at node $i$, $v_i$ represents the quantity of refilling at node $i$, and $waiting_{i}$ represents the waiting time of the sprayer at node $i$. For any node $i\in\mathcal{N}$, we can track the node visited next to it by searching the tuple where node $i$ is the starting node (i.e., $(str_i,end_i)$). Likewise, For any node $i\in\mathcal{N}$, we can track the node visited prior to it by searching the tuple where node $i$ is the end node (i.e., $(str_i,end_i)$). \\

At the beginning of the construction heuristic, $Route_{sp}$ contains only tuples for starting and ending the tour in the depot (i.e., $[\langle (0,0), Q_s,0,0,0 \rangle]$). For given sprayers routing, we can set the value of $l_i$ for each node $i\in\mathcal{N}_f$ simply as $l_i=l_j-q_j+v_j$ if and only if node $i$ is the end node of node's $j$ starting point, namely node $j$ is visited prior to visiting node $i$. However, since we have not established the tender tanker's route yet, we assume that there is a refill at node $j\in\mathcal{N}_f$ if $l_j-q_j < q_i$ for the sequence of nodes $(j,i)$. Therefore, during the construction phase, if $l_j-q_j < q_i$, we set $v_j=Q_s-l_j-q_j$. To calculate the arrival time at node $i\in\mathcal{N}_f$, the value of $y_{i}$ is set to $y_{i} = y_j+s_j+t_{ji}$ if and only if node $i$ is visited following the visit of node $j\in\mathcal{N}$. And finally, we set the value of $m_i$ to 0 for all nodes $i\in\mathcal{N}$. \\

Given the sprayers' routing, we now begin establishing the tender tanker's routing plan. Recall that constructing the sprayers' routes also provides the locations of refill by inspecting nodes where $v_i>0$. The refilling locations are then stored in a list denoted as $\mathcal{RF}$. Nodes within list $\mathcal{RF}$ are the locations to be visited by the tender tanker, to establish the sequence of visits, we order the nodes in an increasing order based on the end of the sprayer's service $s_i+y_i$. Hence, the sequence of nodes to be visited by the tender tanker is now established and the remaining steps are to calculate the arrival time of the tender tanker at each location, the departure time of the tender tanker at each refilling location, and the trips schedule of the tender tanker. The tender tanker refills itself if and only if it's tank level is less than $v_i$ for the next node $i$ in the list $\mathcal{RF}$. The refilling start time and the waiting time at each node are calculated for each node taking into account the synchronization constraints (\ref{eq:FEight7})-(\ref{eq:FTen22}) following the routing of sprayers and the tender tanker.\\

\subsection{Adaptive large neighborhood search}
An extension of \cite{shaw1998using} large neighborhood search (LNS) framework is the adaptive large neighborhood search (ALNS) algorithm. The LNS method is based on constant relaxation and optimization. The approach entails removing certain node visits from the route design. The routing is re-optimized over this relaxed state by reinstalling the node back into the routing plan. If the new re-optimized solution outperforms the current best solution, the former is replaced with the latter and serves as an input for the subsequent iteration.\\

The ALNS framework first introduced by \cite{ropke2006adaptive} and \cite{pisinger2007general} has been used to solve several variants of vehicle routing problems with proven efficient results \cite{azi2014adaptive}, \cite{grangier2016adaptive}, \cite{krebs2021advanced}, and \cite{pfeiffer2022alns}.

ALNS metaheuristic employs the methods of destroy and repair operators; however, contrary to LNS that typically uses a single neighborhood to search, the ALNS works with essentially varied neighborhoods and operators \cite{demir2012adaptive}. Given a set of feasible solutions, the neighborhood is studied by partially destroying the current solution through eliminating certain entries and reinstalling it in each iteration. The deletion and insertion of operators are determined by prior iterations' efficiency. ALNS' adaptive nature adjusts the operator weights based on their performance. Operators with higher performance have a better probability of being selected than those with terrible performance.\\ 

We now describe the ALNS metaheuristic developed in our study. An initial feasible solution $s_0$ is generated using the cluster-first-route-second heuristic approach that we introduced in Section \ref{sec:feasible solution}. We design a set of destroy operators $\mathcal{D}$ and a set of repair operators $\mathcal{Q}$. The destroy operator $d \in \mathcal{D}$ and the repair operator $r \in \mathcal{R}$ are selected probabilistically using the roulette wheel mechanism based on their weights at each iteration. We employ the operators' weight computation method described in the literature \cite{ropke2006adaptive}, \cite{li2016adaptive,alkaabneh2023multiobjective,alkaabneh2023routing}, \cite{ghilas2016adaptive} and \cite{windras2022survey}, see Section \ref{sec:weightUp}. A new neighborhood solution $s_{new}$ is generated in each iteration by using one destroy operator followed by one repair operator. The repaired solution's objective function value is assessed and accepted if it satisfies an acceptance criterion. In particular, if the new neighborhood solution obtained by the destroy and repair operators has an objective function value that is better (i.e., lower for minimization) than the current best solution, we always accept the new solution and update the global best solution. If the new generated solution is better than the current solution and worse than the global best solution, we accept the new solution, make it the current solution, and move forward in the search process of the ALNS (i.e., move to the next iteration). On the other hand, if the new solution is worse than the current solution, it is only accepted with probability following the simulated annealing criteria. Lastly, the algorithm updates the weights of the destroy and repair operators based on the new solution's quality. Contrary to other ALNS algorithms cited in the literature that permit intermediate infeasible solutions to be used to improve their searching capability (e.g., \cite{liu2019adaptive}), our ALNS only accepts feasible solutions.\\

The remaining sections of the ALNS metaheuristic are described as follows: Sections  \ref{sec:Destroy} and  \ref{sec:Repair} present the destroy and repair operators utilized in our ALNS respectively. Section \ref{sec:acceptCriteria} details the acceptance criteria of our ALNS metaheuristic. Section \ref{sec:weightUp} describes the updating of the weights and finally, we present the pseudo code in section \ref{sec:pseudo code}.

\subsubsection{Destroy and repair operators}
\label{sec:Adaptive selection}
A roulette-wheel system controls the selection of the destroy and repair operators of a neighborhood as established earlier. All the destroy and repair operators have an equal chance/probability of being chosen at the initial iteration. After a fixed number of iterations, the probabilities of the destroy operators are updated. 

\subsubsection{Destroy operators}
\label{sec:Destroy}

Our ALNS metaheuristic employs nine distinct destroy operators, each of which implements a unique procedure to fill the removal list $\mathcal{L}$ of nodes to be removed from their positions within each route. A destroy operator accepts a given solution $s$, removes a set of nodes from their positions in their routes, adds these nodes to the removal list $\mathcal{L}$, and returns a solution that is not feasible. The size of list $\mathcal{L}$ is a control parameter that we investigate in Section \ref{ALNS parameter}. The size of the removal list $\mathcal{L}$ is denoted as $p$. We now describe the destroy operators we employ in our ALNS metaheuristic.\\

\textbf{1. Random Removal}\\
The random removal is one of the most basic destroy operators. This operator begins with an empty removal list $\mathcal{L}$. The operator chooses $p$ nodes at random from the selected routes and removes from the solution. The randomization idea gives a diverse search result. The value of parameter $p$ is discussed later.\\

\textbf{2. Waiting Time Removal}\\
This operator removes nodes that have the largest waiting time from the current solution. Nodes are sorted based on the waiting time in an ascending order. The list of eligible nodes to be inserted in the removal list are the nodes within the top 20\% nodes within this sorted list. Finally, $p$ nodes are selected at random to be inserted in the removal list $\mathcal{L}$. \\

\textbf{3. Long-Distance Removal}\\
The operator calculates the traveling time to each node and from that node to the next node (i.e., $t_i=t_{ji}+t_{ie}$ where nodes $j$ and $e$ are the preceding and successive nodes of node $i$, respectively, according to node's $i$ current position). Then the operator sorts the distance in an ascending order and removes $p$ nodes with the longest distance. This saves valuable extra space on the route, which might subsequently take up nodes with shorter distance.\\

\textbf{4. Worst-Order Removal}\\
This destroy operator removes nodes with high-cost savings after removal based on traveling and waiting times at the next node. This operator is specifically designed for VRPMSC-MT by incorporating the waiting time. The cost is defined as the total distance between the node $i$ and its successor (node $e$) and precedence (node $j$) within a specific route. The operator selects and removes node $i = \max_{i\in \mathcal{N}} m_e + t_{ji} + t_{ie}$.\\

\textbf{5. Historical Knowledge Removal}\\
This operator eliminates a node depending on previous information. At every iteration during the execution of the ALNS, the operator calculates the position cost of each node $i$ which is equal to the sum of the distances between its prior and succeeding node and the waiting time at node $i$ and is computed as $t_{i} = t_{ji} + t_{ie}+m_{i}$, where nodes $j$ and $e$ are the preceding and successive nodes of node $i$, respectively, according to node's $i$ current position. The operator updates the node $i$ position cost $t_{i}$ to be the lowest value of all $t_{i}$ values determined up to that iteration. The operator selects a node $j^{*}$ on a route with the maximum difference from its node cost i.e. $j^{*} = max_{j\in\mathcal{N}} {t_{j} - t_{j^*}}$ and adds to the removal list $\mathcal{L}$, the procedure continues until $p$ nodes have been added to the list.\\

\textbf{6. Route Removal}\\
The operator accepts the set of sprayers' routes as input and randomly removes all nodes within a chosen route from a set of routes in the solution into a removal list $\mathcal{L}$.\\

\textbf{7. Zone Removal}\\
The zone removal operator eliminates a collection of nodes from given routes based on their placement in the map. This operator generates a point in the map at random and adds all nodes within a predetermined radius of the selected point to the removal list $\mathcal{L}$. In our work, we set the value of the radius to be 0.25 miles. \\

\textbf{8. Proximity Removal}\\
A subclass of the Shaw removal operator is proximity removal. The Shaw removal operator eliminates a set number of relatively predefined similar nodes because reshuffling similar nodes increases the probability of generating feasible and improved solutions. The core idea of Shaw removal was introduced by \cite{shaw1998using}. This proximity removal operator picks a node at random and adds the selected node and four nodes related to the selected node to be added to the removal list. The four selected nodes are the two preceding and the two following the selected node within its route. The number of times we run this destroy operator is set to be $5.0\%\vert \mathcal{N}_f\vert$.\\

\textbf{9. Refill Position Removal}\\
This method is specially designed for our VRPMSC-MT. It adds all refilling locations to the removal list $\mathcal{L}$, namely, all nodes in list $\mathcal{RF}$ are added to the removal list $\mathcal{L}$.

\subsubsection{Repair Methods}
\label{sec:Repair}
In this section, we present the two main greedy repair operators we use in our ALNS metaheuristic. Repair operators are used to repair the partially destroyed solution by re-inserting the nodes in the removal list $\mathcal{L}$ sequentially into positions within sprayers. The two repair operators we used were inspired by \cite{ropke2006adaptive}. See the brief description below:\\

\textbf{1. Greedy Repair:}\\
This repair operator receives the removal list $\mathcal{L}$ and iteratively selects a node from the list to put one at a time into existing routes in the best feasible position. The cost of each insertion is calculated at each iteration to obtain the insertion with the minimum cost. The insertion position with the lowest insertion cost is chosen, and the corresponding insertion operation is performed.\\

\textbf{2. Regret Repair}\\
The regret repair operator calculates for each repair request a regret value equal to the objective function difference between the best insertion position and the second best insertion position. As a result, repair requests with a high level of regret will be inserted first i.e. the basic regret technique chooses the insertion that would be most regretted if it is not inserted now. After every repair, the regrets must be recalculated because repaired nodes are no longer available.\\

During the execution of the ALNS metaheuristic, any time new routes for sprayers are generated via the repair operation, we use the mechanism described in Section \ref{sec:soldRep} to generate refilling locations and the trips of the tender tanker to obtain a complete feasible solution.

\subsubsection{Acceptance Criteria}
\label{sec:acceptCriteria}
Denote a new feasible solution by $s_{new}$ that was obtained by implementing a destroy operator and repair operator on a given solution, the global best solution found so far as $s^*$, and the current solution as $s_{current}$, and the objective function value of a feasible solution as $f(\,)$. The new solution $s_{new}$ is always accepted if $f(s_{new})<f(s^*)$. The new solution is also accepted if $f(s_{new}) < f(s_{current})$. On the other hand, the new solution might be accepted if $f(s_{new}) > f(s_{current})$ using a probability function defined as $e^{(f(s_{new}) - f(s_{current})/Temp}$ and controlled by a parameter of $Temp$, where $Temp$ is the temperature. After each iteration, the value of $Temp$ reduces following the formula $Temp = Cooling*Temp$, where $(0 < Cooling < 1)$ represents the cooling rate. The formula $Temp_{0} = \rho * f(s_{0})$ is used to compute the initial $Temp$ where $s_0$ is the initial feasible solution obtained from Phase 1 of the metaheuristic (see Section \ref{sec:feasible solution}). Note that if there is no improvement to the global best solution after a certain number of iterations, the use $s^*$ as the current solution and the ALNS moves forward. Lastly, our ALNS terminates after a certain number of iterations. 

\subsubsection{Intensive Local Search}
\label{sec:localSearch}
In our ALNS to further intensify the search process we use a mathematical model to improve on some solutions that are promising. More specifically, whenever a new solution is better than the current best solution, we explore the neighbor of that solution to further improve that solution by solving mathematical model (\ref{eq:FObjt})-(\ref{eq:FIntg}) with certain binary variables being fixed. Note that finding optimal routes of the sprayers alone does not guarantee finding the optimal solution to model (\ref{eq:FObjt})-(\ref{eq:FIntg}) since finding the optimal refilling nodes and the route of the tender tanker is an optimization problem by itself. That said, each time during the execution of the ALNS a new best solution is found, we fix the sprayers' routing variables (i.e., $x_{ij}$) using the values found in the new best solution and we solve model (\ref{eq:FObjt})-(\ref{eq:FIntg}) with fixed variables. This strategy will guarantee finding the optimal refilling nodes and the routes of the tender tanker given feasible sprayers' routes. However, since this MIP is also NP-Hard we only run the solver with a time limit of 60 seconds. 

\subsection{Updating the weights}
\label{sec:weightUp}

The ALNS algorithm updates the weights of the operators based on the performance of each operator during the search process. In our work, the weight of each operator is updated at the end of a fixed number of iterations (called a segment) rather than at the end of each iteration. Note that since we only use two repair operators, the weight of each repair operator does not update, rather we only update the weight of the destroy operators. \\

Based on the performance of a destroy operator $d\in\mathcal{D}$ during the search process in any segment, a score $\eta_d$ is assigned to it at the end of that segment. The formula for updating the value of $\eta_d$ is as follows:
$$\eta_o^{seg-1} = \lambda\eta_o^{seg-1}+ (1-\lambda)\psi$$
where $\lambda\in(0,1)$ is a smoothing parameter and $\psi$ is the performance score of a destroy $d\in\mathcal{D}$ operator and it is assigned as:

\begin{equation*}
  \psi = \left \{
  \begin{aligned}
    &\psi_1, && \text{if the new solution is better than the best solution found so far}\\
    &\psi_2, && \text{if the new solution is better than the current solution}\\
    &\psi_3, && \text{if the new solution is accepted}\\
    &\psi_4, && \text{if the new solution is rejected}
  \end{aligned} \right.
\end{equation*}
with $\psi_1\geq \psi_2\geq \psi_3\geq \psi_4$, in our implementation we use $\psi_1=7, \psi_2=4, \psi_3=2,$ and  $\psi_4 = 1$. The initial weight $\eta_d$ at the first iteration is set to 1 for all operators. At the end of each segment after calculating the values of $\eta$, the probability assigned to each operation $d$ is $ps_d = \frac{\eta_d}{\sum_{d\in \mathcal{D}}\eta_d}$

\subsubsection{ALNS Framework: pseudo code}
\label{sec:pseudo code}

\begin{algorithm}[H]
  \renewcommand{\arraystretch}{0.5}
  \small
\DontPrintSemicolon
  \KwInput{A feasible solution $s_0$ produced by a construction heuristic and a set of ALNS parameters: $Temp,\ \mathbf{\psi},\ iterMax,\ MaxNoImprov$}
  \KwOutput{The best known solution $s^*$}
  Set $iter \gets 0,\ s^* \gets s,\ s_{current} \gets s $, $\eta = [1,...,1]$\;
  \While{$iter < iterMax$}
   {
   		Select a destroy operator $d$ and repair operator $r$ by a roulette wheel mechanism using weights $\mathbf{ps}$\;
   		Apply destroy operator $d$ on $s_{current}$ and get $s^\prime$\;
   		Apply repair operator $r$ on $s^\prime$ and get $s_{new}$\;
   		\If{$f(s_{new}) < f(s^*)$}
    {
        $s^* \gets s_{new},\ s_{current} \gets s_{new}$\;
        \label{line8}
        Implement intensive local search and obtain $s^*$ based on the solution of the mathematical model\; \tcp*{See Section \ref{sec:localSearch}}
        $NoImproveCounter \gets 0$\;
        
    }
    \Else
    {     $NoImproveCounter \gets NoImproveCounter + 1$\;
    	  \If{$f(s_{new}) < f(s_{current})$}
    {
        $s \gets s_{new}, s_{current} \gets s_{new}$\;
        \label{line13}
        
    }
    \Else
    {
    	$probAccept \gets e^{\frac{f(s_{new})-f(s_{current})}{Temp}}$\;
    	Generate a random number $U\sim Uniform[0,1]$\;
    	\If{$U < probAccept$}
            {
                $s_{current} \gets s_{new}$    \tcp*{Accept the new solution which is worst than the current solution for diversification}

            }
            \Else
            {
            	$s_{current} \gets s_{current}$ \tcp*{Reject the new solution and use the current solution}
            	
            }
    }
    }
    $iter\gets iter+1$\;
    \If{$iter \% segmentLength == 0$}
        {Update $\mathbf{\eta}$ and calculate $\mathbf{ps}$\;}

    \If{$NoImproveCounter \% MaxNoImprov == 0$}
        {$s \gets s^*,\ NoImproveCounter \gets 0$}
   }
\caption{ALNS framework for the VRPMSC-MT}\label{ALNSalgo}
\end{algorithm}

STSARP
\section{Computational Results}
\label{sec:Computational}
The goals of computational results and analysis are three: (1) we first focus on presenting numerical analyses to demonstrate the effectiveness of our proposed ALNS meta-heuristic combined with local search, (2) to compare and contrast different models presented in Section \ref{sec:mathModel}, and (3) to demonstrate the benefits of using an optimization approach to solve the sprayer-tanker routing problem as opposed to the approach implemented in practice. Given the fact that our work is in collaboration with a private company, the instances we generate follow random number generators, with specifics presented in each of the sections discussed below.\\

The model and ALNS met-heuristic presented in this paper are coded in Python and we use Gurobi 9.1.1 optimization software as the MILP solver. All the experiments were conducted on a computer with Intel(R) Core(TM) i7-8700 CPU @ 3.20GHz processor with 16.0GB RAM.\\

In farm spraying, the quantity of fertilizer to spray each node can be calculated in advance based on several factors such as crop yield and effective pest control. On the other hand the consequences of excessive chemical application results in damages such as major soil degradation, nitrogen leaching, soil compaction, reduction in soil organic matter, and loss of soil carbon.

Regarding sprayer type, for example, the boom sprayers have $25 - 100$ gallon tank capacity, $2 - 5$ gallon per minute nozzle flow rates and maximum pressures of 60 psi. In our computational analysis, we use the parameters of the boom sprayer and we assume that the farmer uses a boom sprayer. To generate instance parameters, we truncate some of the values to maintain the privacy of data, we set the sprayer maximum tank capacity ($Q_s$) to 25 units of fertilizer and the demand of fertilizer $q_i$ for each node to be an integer number with the range of $2 - 5$ units of fertilizer.

In all our experimental analysis we randomly generated various test instances for a wide range of situations since there are no existing benchmark examples in the literature that are relevant to our problem setting because the sprayer-tanker routing problem with multiple trips is a new variant of VRP. Our experimental test-bed consists of the following parameters:

\begin{itemize}
\item Number of nodes to be sprayed $\vert \mathcal{N}_f\vert$: 15, 16, 17, ...,20 for small instances, 25, 30, 35, and 40 for medium instances, and 50 and 60 nodes for large-scale instances. For big farms, the number of nodes to be sprayed is less than 60 nodes. Hence, our generated intances are realistic.
\item Number of sprayers: $\lbrace 2,3,4,5\rbrace$.
\item Tender tanker capacity $Q_t$: 100 units of fertilizer.
\item Service time of Sprayer $s_i$ in units of time is consistent with the demand for fertilizer of each node.
\end{itemize}

\subsection{ALNS parameter setting}
\label{ALNS parameter}
In this section, we perform extensive analysis to find good values for ALNS parameters to use. Our analysis focuses on some parameters that affect the performance of the ALNS. These parameters are (i) the size of the removal set $\mathcal{L}$ (i.e., how many nodes to remove when performing a destroy operator), (ii) the number of ALNS iterations before termination, and (iii) which destroy operators to include. Once we identify the best value for these parameters, we fix these values.  

\subsubsection{Size of removal list $\mathcal{L}$}
In this section, we assess the efficiency and effectiveness of varying the size of the removal list $\mathcal{L}$ to decide on the number of nodes to be removed during a destroy operator. The number of nodes selected to be removed by the chosen removal operator in each iteration is an important factor affecting the performance of the ALNS. If only a few nodes are removed from the current solution, it may be impossible to avoid a given neighborhood local optimum. On the other hand, removing a significant number of nodes makes it expensive to get a feasible solution in terms of computational time, resulting in longer computational time. The ALNS metaheuristic is designed to have rapid convergence, although diversification is critical at the initiation of the ALNS algorithm to ensure a large range of solution search space. To this end, we investigate three unique sizes of the removal list, the first size of removal list at any ALNS iteration will be a random integer number between $5\%$ and $10 \%$ of the total number of nodes denoted as RL1, the second size is $7\% - 15 \%$ of number of nodes denoted as RL2, and finally the third size is $7.5\%-12.5\%$ denoted as RL3.

We perform experiments to determine the best size of the removal list for tuning our ALNS metaheuristic. We generate 12 random instances for various number of nodes and number of sprayers. For each instance, we generate 15 runs resulting in a total of 180 runs. Other model parameters remain constant. The average computational time and percentage gap for each removal list size are estimated and compared. The detailed results of the different sizes of removal list of the ALNS metaheuristic for different test instances and number of sprayers are shown in Table \ref{Bounds}. We label an instance by encoding them as Instance-number of nodes-number of sprayers. For instance, I-20-S2 means an instance with 20 nodes and two sprayers. The table format is as follows: column one shows the instance label and the remaining columns display the results of the computational time (in seconds) and percentage improvement (Gap \%). The gap is calculated in reference to the improvement in the objective function value of the best solution returned by the ALNS against the initial feasible solution. Note that each value in the table is calculated based on the average of 15 runs for each instance.  

The results show that the ALNS computational time increases as the size of the removal list increases for a given instance. This result is consistent with the intuition that when a larger set of nodes is removed, more computational time is needed to repair the destroyed solution. For example, the results of the test instances I-20-S2 show an average computation time of 22.36, 47.49, and 70.71 seconds for RL1, RL2, and RL3 respectively. More essentially, we notice that for all the test instances the computational time for RL1 (i.e., the size of the removal list is between 5\% - 10 \% of the number of nodes) was much lower.

On the other hand, the results regarding the gap are not as clear as the results regarding the computational time. More specifically, RL1 does not always find the highest quality solutions. In some cases, we see that RL2 achieves the best results (for instance, in the case of I-30-S2). Likewise, the percentage gap is higher for test instances I-20-S4, I-30-S4, I-40-S3 and I-40-S4 for RL1, I-20-S2, I-20-S3, I-20-S3 and I-40-S2 for RL2 and I-25-S2, I-25-S3, I-25-S4 and I-30-S3 for RL3. However, the overall average indicates that RL1 has a better percentage improvement compared to RL2 and RL3. The results show that higher computational time does not necessarily lead to improved solutions.

\begin{figure}[!ht]
\centering
\includegraphics[width=14cm]{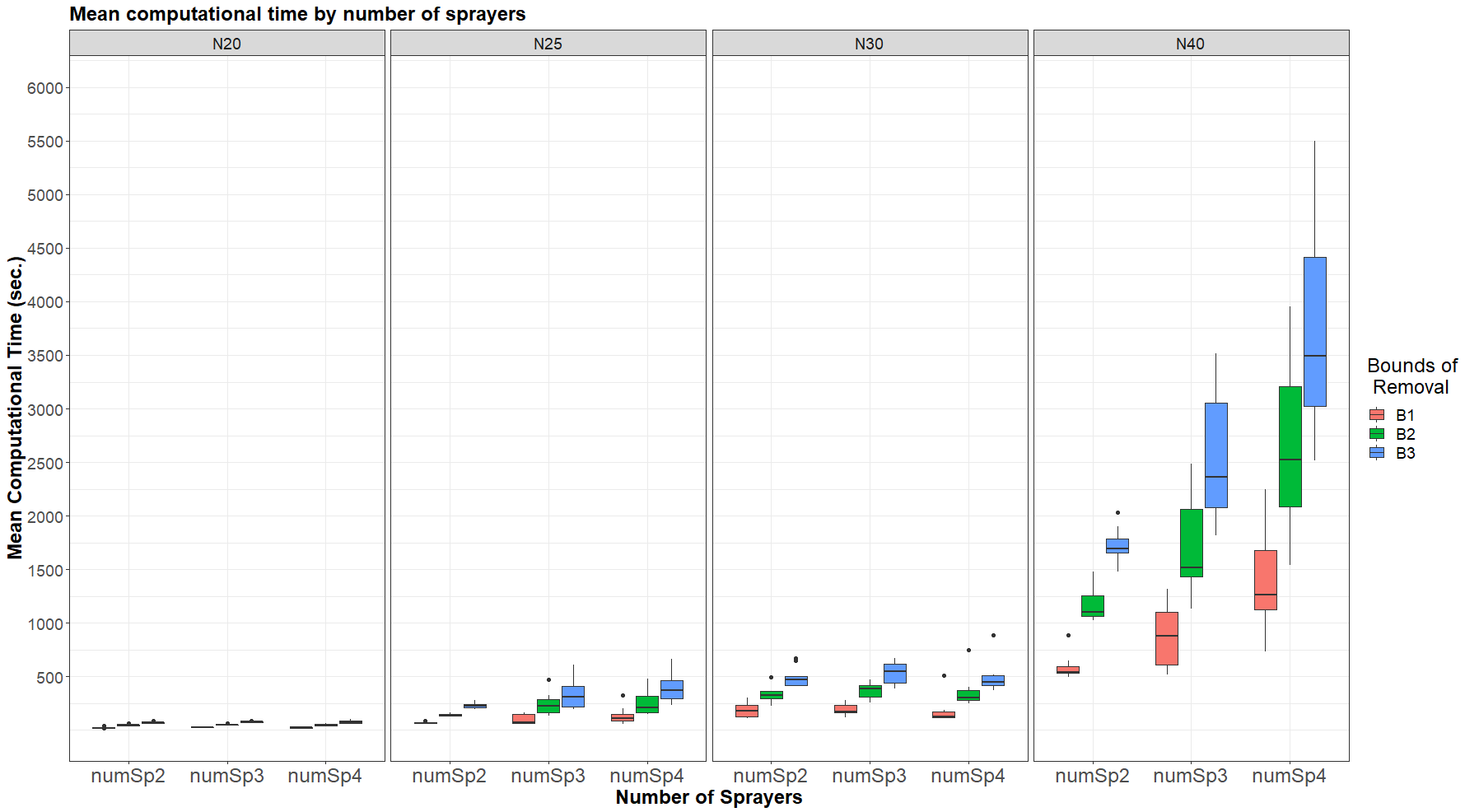}
\caption{Results on mean ALNS computational time for different instance and sprayers for the Bounds of removal}
\label{fig:ALNS Time}
\end{figure}

\begin{figure}[!ht]
\centering
\includegraphics[width=14cm]{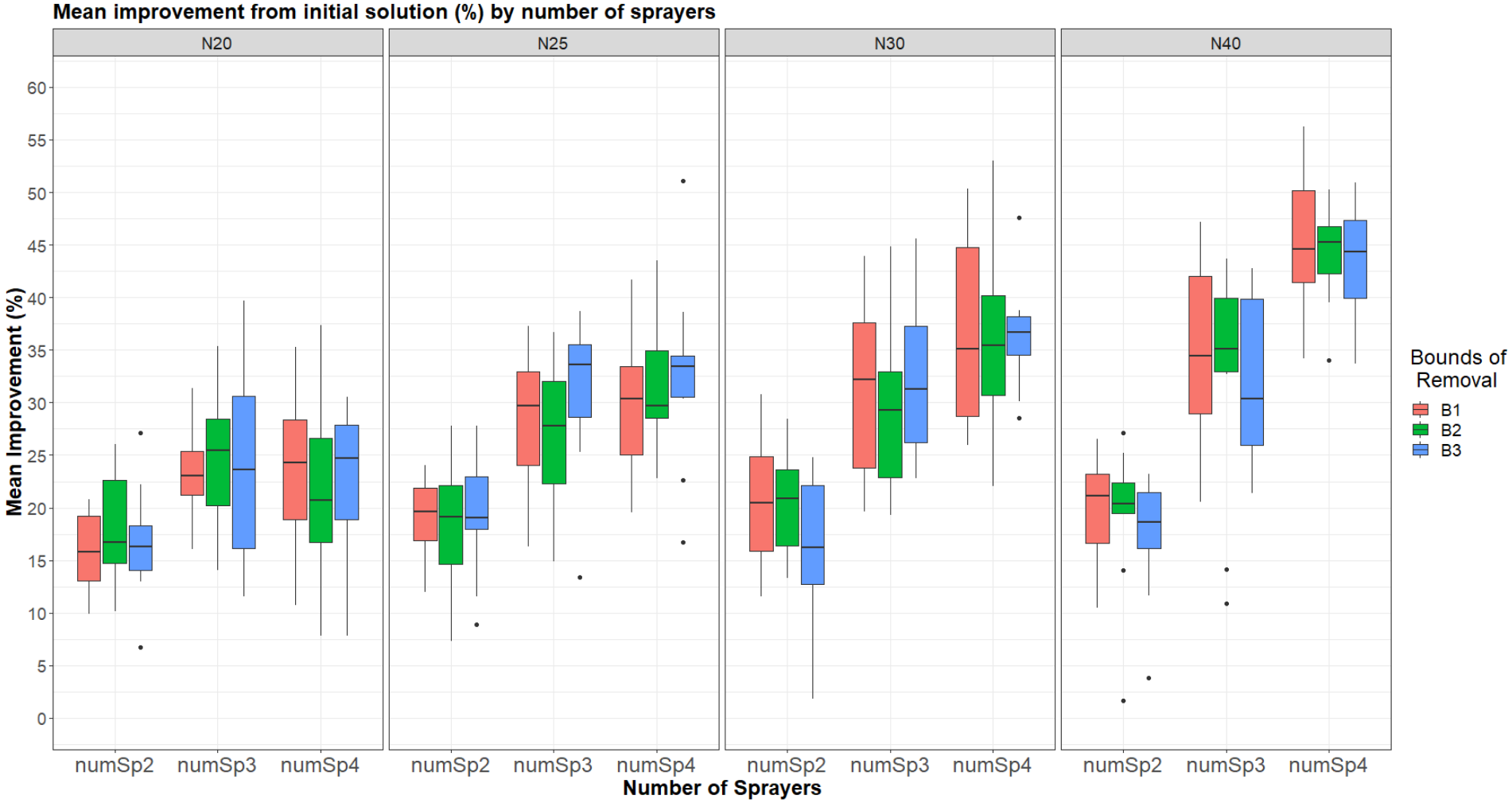}
\caption{Results of the average percent ALNS improvement for different instances and sprayers for the Bounds of removal}
\label{fig:Improvement}
\end{figure}

To present the complete results summarized in Table \ref{Bounds}, Figures \ref{fig:ALNS Time} and \ref{fig:Improvement} are created to show the box-plots of the gap and computational time from the data points presented in Table \ref{Bounds}. More specifically, Figure \ref{fig:ALNS Time} illustrates the ALNS computational time for different instances presented in Table \ref{Bounds}. We note that RL1 yields the lowest computational time and the RL3 yields the highest. On the other hand, Figure \ref{fig:Improvement} shows the percentage of ALNS improvement for different instances and the number of sprayers. We observe that the performance of each of the sizes of the removal list was relatively different across all the instances and the number of sprayers. Since RL1 gives better results in terms of the percentage improvement and computational time, it is reasonable to choose that range for all of the computational analysis we perform. Therefore, we opt to use a removal list of size $RL 1:5\% - 10 \%*\vert \mathcal{N}_f\vert$ as the default in order to maintain a balance for computational time and solution quality.

\begin{table}[]
\begin{center}
\caption{Results of various removal list size with different instances.}\label{Bounds}
\begin{tabular}{llllllll}
& \multicolumn{3}{c}{Gap(\%)}&  & \multicolumn{3}{c}{Time (sec)}\\ 
\cline{2-4} \cline{6-8} 
Inst.& RL1 & RL2 & RL3 &  & RL1 & RL2 & RL3 \\ 
\hline
I-20-S2 & 15.74 & \textbf{18.01} & 16.65&  & \textbf{22.36}   & 47.49   & 70.71 \\
I-20-S3 & 23.79 & \textbf{24.59} & 24.29&  & \textbf{24.20}   & 49.85   & 76.80   \\
I-20-S4 & \textbf{23.87} & 22.45 & 22.49&  & \textbf{22.92}   & 51.40   & 78.83   \\
I-25-S2 & 18.94 & 18.51  & \textbf{19.11} &  & \textbf{69.12}   & 142.51  & 226.60  \\
I-25-S3 & 28.40 & 26.66 & \textbf{30.95} &  & \textbf{99.84}   & 240.33  & 332.09  \\
I-25-S4 & 29.74& 31.59  & \textbf{32.59} &  & \textbf{133.49}  & 254.63  & 387.34  \\
I-30-S2 & 20.37  & \textbf{20.44} & 16.01 &  & \textbf{190.55}  & 346.19  & 491.74  \\
I-30-S3 & 31.29 & 29.38  & \textbf{32.15} &  & \textbf{190.59}  & 370.65  & 535.07  \\
I-30-S4 & \textbf{36.58} & 35.98& 36.37 &  & \textbf{175.97}  & 352.81  & 490.51  \\
I-40-S2 & 19.89& \textbf{19.20} & 17.27 &  & \textbf{586.94}  & 1165.77 & 1720.92 \\
I-40-S3 & \textbf{34.82} & 32.42& 32.01 &  & \textbf{873.63}  & 1725.63 & 2543.18 \\
I-40-S4 & \textbf{45.38} & 43.99& 43.34 &  & \textbf{1379.02} & 2596.77 & 3742.19 \\ \hline
Avg.    & 27.40 & 26.93 & 26.94  &  & 314.05& 612.00  & 891.33  \\
Min     & 15.74 & 18.01 & 16.01 &  & 22.36  & 47.49   & 70.71   \\ 
Max.    & 45.38 & 43.99 & 43.34 &  & 1379.02& 2596.77 & 3742.19 \\

\hline
\end{tabular}
\end{center}
\end{table}

\subsubsection{Maximum number of iterations before termination}
We also investigate the maximum number of iterations the ALNS performs before stopping as the termination criteria of our ALNS depends on the maximum number of iterations. We perform the experiments with four distinct maximum number of iterations: $50*\vert \mathcal{N}_f\vert,\, 100*\vert \mathcal{N}_f\vert,\, 150*\vert \mathcal{N}_f\vert,$ and $200*\vert \mathcal{N}_f\vert$.

We generate random instances for varying numbers of nodes and number of sprayers, and for each instance, we generate 15 runs, resulting in a total of 135 runs. The average of computational times and the average percentage improvement are computed. Table \ref{maxIterTable} summarizes the results of comparing the maximum number of iterations of the ALNS meta-heuristic before termination, Table \ref{maxIterTable} has the same format and structure as Table \ref{Bounds}.   

\begin{table}[]
\caption{Results of various maximum number of iterations with different instances.} \label{maxIterTable}
\begin{tabular}{llllllllll}
& \multicolumn{4}{c}{Gap (\%)}  &  & \multicolumn{4}{c}{Time (sec)}\\ \cline{2-5} \cline{7-10} 
Inst.& $50*n$  & $100*n$ & $150*n$& $200*n$ &  & $50*n$  & $100*n$ & $150*n$ & $200*n$   \\ \hline
I-30-S2 & 25.35 & 25.84 & \textbf{25.90} & 25.87 &  & \textbf{189.67}  & 479.56  & 816.88  & 1201.85 \\
I-30-S3 & 30.29 & 30.69 & \textbf{30.71} & 30.79 &  & \textbf{177.65}  & 433.73  & 741.87  & 1117.77 \\
I-30-S4 & 25.21 & 25.74 & 25.90 & \textbf{25.93} &  & \textbf{173.03}  & 447.04  & 802.78  & 1161.33 \\
I-40-S2 & 20.31 & 20.53 & \textbf{20.57} & 20.53 &  & \textbf{419.03}  & 1043.39 & 1862.72 & 2865.29 \\
I-40-S3 & 33.34 & 33.44 & 33.60 & \textbf{33.66} &  & \textbf{679.38}  & 1499.19 & 2542.08 & 3759.90 \\
I-40-S4 & 45.49 & \textbf{46.33} & 46.15 & 46.23 &  & \textbf{833.19}  & 1990.26 & 3277.26 & 4933.31 \\
I-50-S2 & 24.25 & 24.31 & \textbf{24.45} & 24.44 &  & \textbf{905.55}  & 2310.43 & 4194.76 & 6605.42 \\
I-50-S3 & 40.87 & 41.38 & 41.25 & \textbf{41.39} &  & \textbf{1278.06} & 2976.36 & 5140.58 & 7861.73 \\
I-50-S4 & 46.04 & 46.51 & \textbf{46.71} & 46.38 &  & \textbf{1768.30} & 4088.58 & 6813.81 & 9842.66 \\ \hline
Avg.    & 32.35 & 32.75 & 32.80 & 32.80 &  & 713.76  & 1696.51 & 2910.30 & 4372.14 \\
Min     & 20.31 & 20.53 & 20.57 & 20.53 &  & 173.03  & 433.73  & 741.87  & 1117.77 \\
Max     & 46.04 & 46.51 & 46.71 & 46.38 &  & 1768.30 & 4088.58 & 6813.81 & 9842.66 \\ \hline
\end{tabular}
\end{table}

\begin{figure}[!ht]
\centering
\includegraphics[width=14cm]{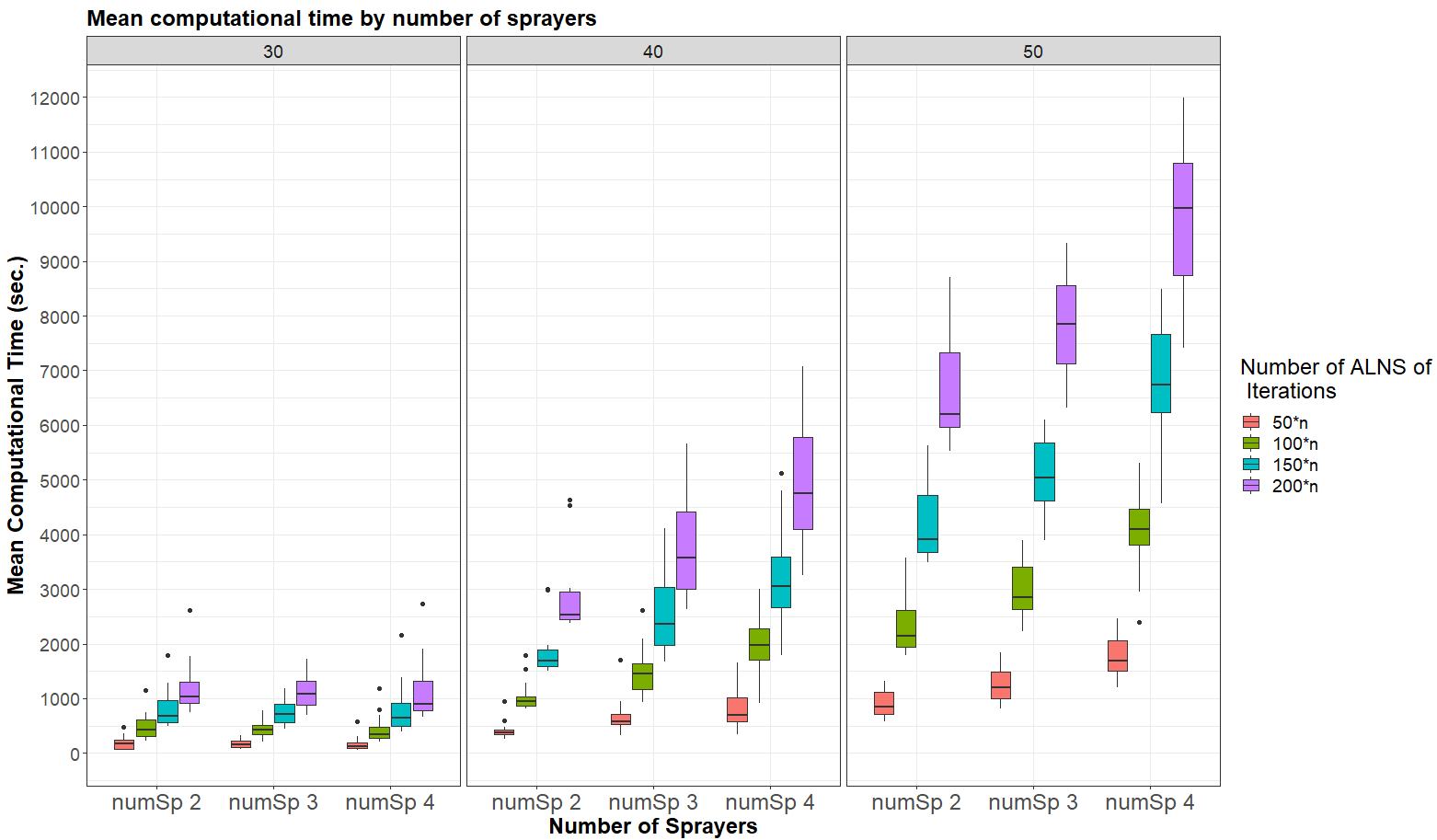}
\caption{Results on average ALNS time for different instances and sprayers for different number of ALNS iterations}
\label{fig:IterMaxTime}
\end{figure}

\begin{figure}[!ht]
\centering
\includegraphics[width=14cm]{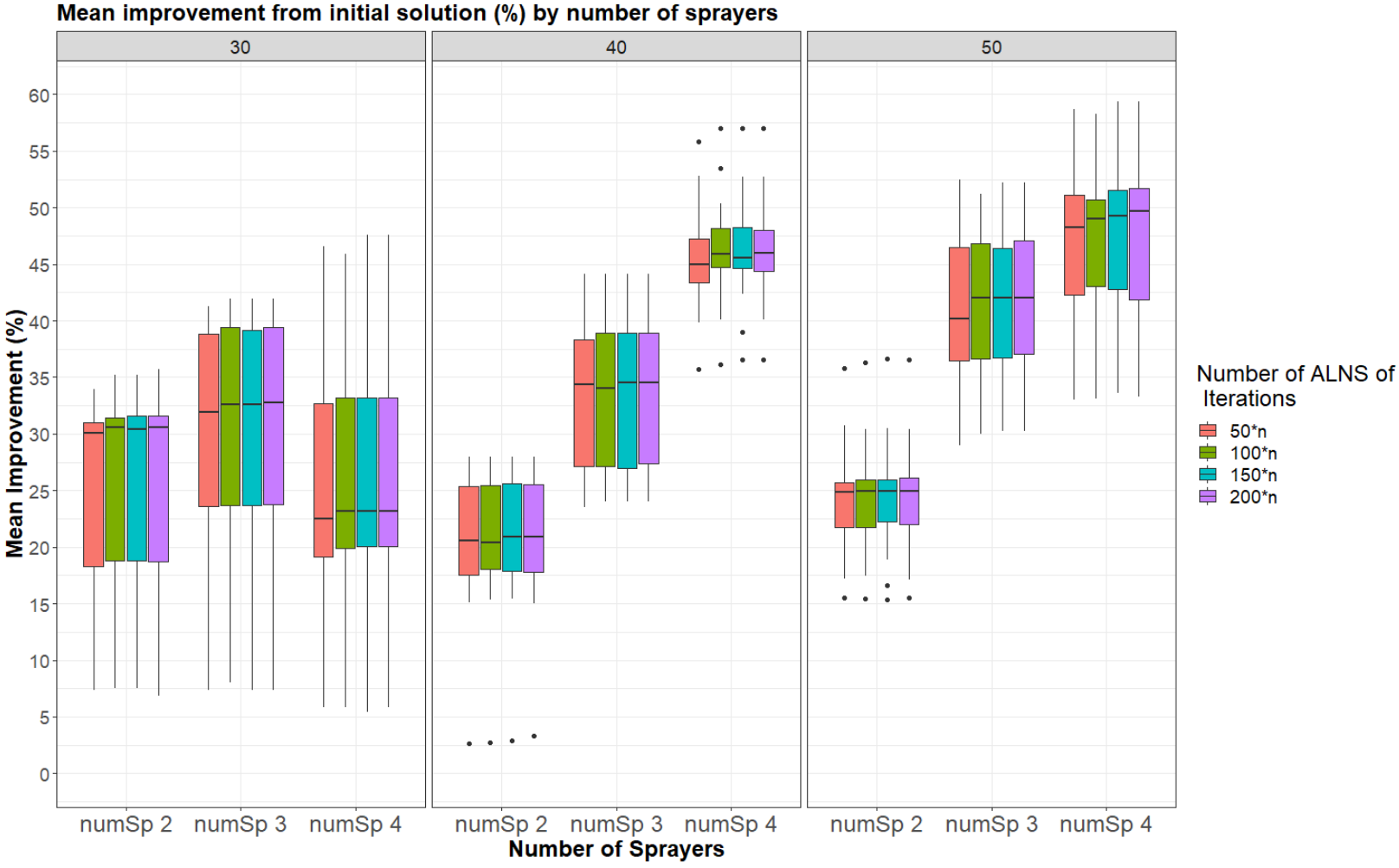}
\caption{Results on average percent ALNS improvement for different instance and sprayers for different ALNS iterations}
\label{fig:Maxpercent}
\end{figure}

Figure \ref{fig:IterMaxTime} illustrates the average computational time for different number of sprayers and instances across all four types of iterations. The $numSp$ indicates the number of sprayers and the values at the top horizontal panel show the number of nodes. The results show that for the same number of nodes and the same number of sprayers, the computation time increases drastically with increasing number of maximum number of iterations. Figure \ref{fig:Maxpercent} displays the results for the average percent improvement for the different iterations. We observe that for the same instance and the same number of sprayers, we do not notice any substantial difference between the percent mean improvement from the initial solution in all four iterations types. For example with 30 nodes, the percent mean improvement for iterations $50*\vert \mathcal{N}_f\vert,\, 100*\vert \mathcal{N}_f\vert,\, 150*\vert \mathcal{N}_f\vert$ and $200*\vert \mathcal{N}_f\vert$ are approximately the same across all the different number of sprayers. We observe a similar pattern in the test with 40 nodes and 50 nodes as well. Based on the results, it is more reasonable for to set the maximum number of iterations to be $100*\vert \mathcal{N}_f\vert$ for our computational analysis as it provides a good effectiveness considering the computational time.

\subsubsection{Destroy Operators}
In this section of the paper, we investigate the effectiveness of all destroy operators to obtain insight into which destroy operators contribute most to the overall optimal solution generation.

We generate 9 random instances with different number of nodes and number of sprayers. The experimental run for each instance follows the same method in section \ref{ALNS parameter}. We set the ALNS parameters as follows: removal list size is RL1: $5\% -10 \%$, maximum iteration of $100*\vert\mathcal{N}_f \vert$. Other model parameters remain unchanged. 

We measure the performance of a destroy operator based on the percentage of iterations in which the repaired solution was better than the current best solution divided by the number of times the destroy operator was selected during the ALNS execution. This value helps to determine which operators are most likely to obtain a new best solution. More specifically, if we denote that performance quantity by $Performance$, then $Performance = NewBestCounter/UsageCounter*100\%$ where $NewBestCounter$ is the number of times the new repaired solution is better than the current best solution and $UsageCounter$ is the number of times the removal operator was selected. 

Intuitively, the most frequently used operators are the ones that contribute the most to obtaining a new best solution; however, the most frequently used operators are not always the ones that obtain the new best solution. We present the results in Table \ref{Operators}. The format of Table \ref{Operators} is as follows: the first column shows the instance label as we describe in section \ref{ALNS parameter} and the subsequent set of columns show the different destroy operators we use. We also compute the overall average performance $(Avg)$, the minimum $(min)$ as well as the maximum $(max)$ of each removal operator.

We notice that the route removal, historical knowledge removal, and worst distance removal operators contribute the most to the new best solutions, with averages of 1.44\%, 1.22\%, and 1.11\%, respectively, while waiting time removal, random removal, and proximity removal operators contribute more to solution diversity with averages of 0.17\%, 0.43\% and 0.73\% respectively. 

\begin{landscape}
\begin{table}[]
\begin{center}
\caption{Results on comparative performance of removal operators. \label{Operators}}
\begin{tabular}{llllllllll}
&\multicolumn{9}{c}{Removal Operators}\\ \cline{2-10} 
Inst. & Random & \begin{tabular}[c]{@{}l@{}}Waiting\\ Time\end{tabular} & \begin{tabular}[c]{@{}l@{}}Longest\\ Distance\end{tabular} & \begin{tabular}[c]{@{}l@{}}Worst\\ Distance\end{tabular} & \begin{tabular}[c]{@{}l@{}}Historical\\ Knowledge\end{tabular} & Route & \begin{tabular}[c]{@{}l@{}}Refill\\ Position\end{tabular} & Zone & Proximity \\\hline
I-30-S2 & 0.23& 0.07  & 0.70 & 0.73 & 0.55 & 2.01 & 0.55 & 0.50 & 0.36 \\
I-30-S3 & 0.44& 0.17 & 0.60 & 1.30 & 1.17 & 1.46 & 1.08 & 1.07 & 0.93    \\
I-30-S4 & 0.45& 0.19 & 0.90  & 1.61 & 1.48 & 1.38 & 0.48 & 0.99 & 1.17   \\
I-40-S2 & 0.17& 0.08 & 0.39 & 0.43 & 0.93 & 1.17 & 0.47 & 0.68 & 0.40    \\
I-40-S3 & 0.43& 0.09 & 0.77 & 1.38 & 1.14 & 1.46 & 0.89 & 0.76 & 0.62    \\
I-40-S4 & 0.72& 0.40 & 1.27 & 1.63 & 1.63 & 1.15 & 1.03 & 1.26 & 0.93    \\
I-50-S2 & 0.25& 0.02 & 0.65 & 0.73 & 1.04 & 1.23 & 0.50 & 0.61 & 0.48    \\
I-50-S3 & 0.48& 0.14 & 0.87 & 0.94 & 1.33 & 1.36 & 0.95 & 0.81 & 0.68    \\
I-50-S4 & 0.72& 0.38 & 1.07 & 1.24 & 1.75 & 1.76 & 1.40 & 1.09 & 1.02    \\ \hline
Avg. & 0.43& 0.17 & 0.80 & 1.11 & 1.22 & 1.44& 0.82 & 0.87 & 0.73    \\
Min.& 0.17& 0.02 & 0.39 & 0.43 & 0.55 & 1.15& 0.47 & 0.50 & 0.36     \\
Max.& 0.72& 0.40 & 1.27 & 1.63 & 1.75 & 2.01 & 1.40 & 1.26 & 1.17    \\ \hline
\end{tabular}
\end{center}
\end{table}
\end{landscape}

\begin{figure}[!ht]
\centering
\includegraphics[width=10cm]{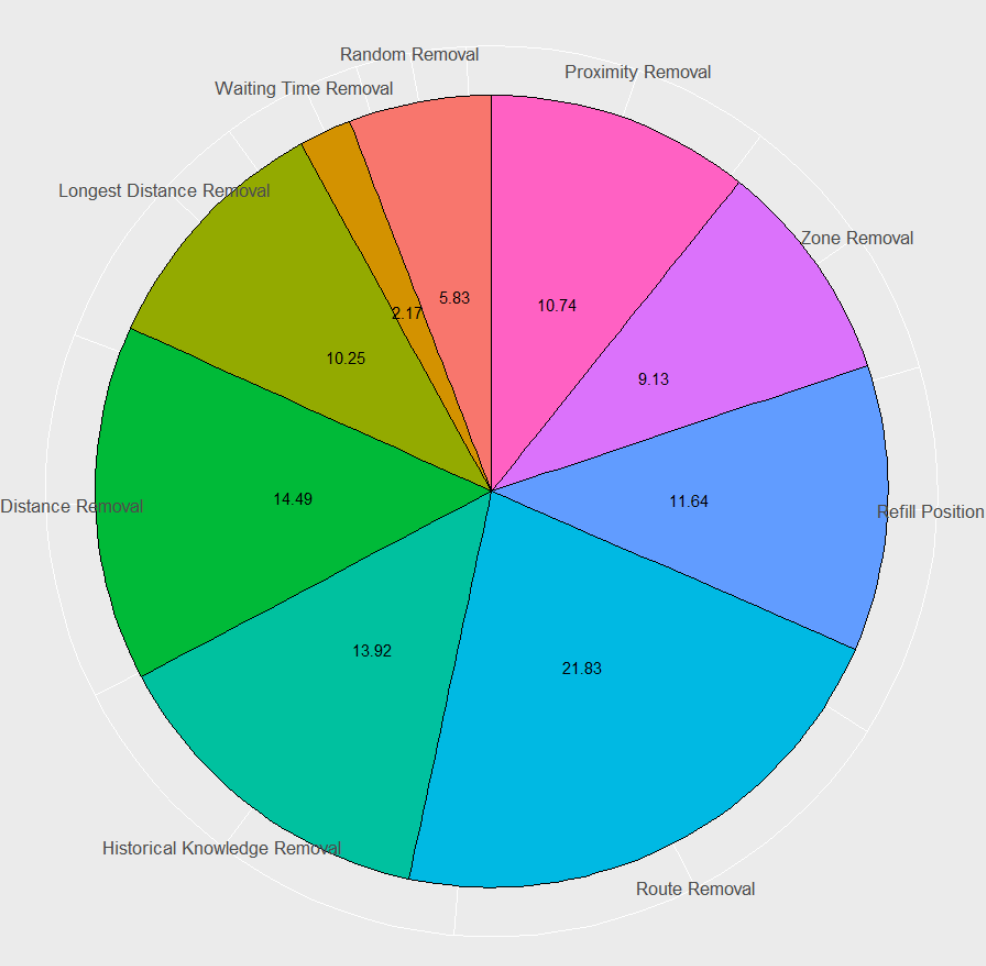}
\caption{Percentage contribution of each removal operator in getting a new solution that is better than the current best solution.}
\label{fig:removalP}
\end{figure}

Normalizing the data presented in Table \ref{Operators} we obtain the percentage contribution of each destroy operator in getting a new solution that is better than the current best solution, the results are displayed in Figure \ref{fig:removalP}.

\subsection{Computational Performance of the Metaheuristic}

In this section, we illustrate the computational performance of the developed metaheuristic in terms of computational time, objective function value, optimality gap, and percentage improvment of the ALNS by comparing the best solution found during the search against the initial feasible solution used as an input. We generate 16 instances at random, with different number of nodes and number of sprayers. For each instance, we generate 10 experimental runs. When assessing the performance of our ALNS metaheuristic, we compare the meta-heuristic against the Gurobi solver by solving the mathematical model (\ref{eq:FObjt})-(\ref{eq:FIntg}). We present the results in Table \ref{SolverALNS}. The first column shows the instance label, the subsequent columns 2-4  show the computational time in seconds, objective function value $Obj$, $Gap\%$ for Gurobi solver, and columns 5 - 7 show the computational time in seconds, objective function value, the percentage ALNS improvement from the initial solution $(Imp\%)$ corresponding to the ALNS metaheuristic. For all instances, the time limit for solving the mathematical model using Gurobi is set to 7,200 seconds.

The computational analysis we present in Table \ref{SolverALNS} indicates an excellent performance of the developed ALNS metaheuristic in finding high quality solutions in a shorter computational time when compared to the solver. On average, the metaheuristic yields better solutions than the solver in proportion to the time required by the solver. Specifically, the metaheuristic finds such high-quality solutions in a much shorter computing period. For example, with 15 nodes and 2 sprayers, the solver requires 3,888.86 seconds to reach an optimal objective function value of 48.19  but the metaheuristic requires only 6.78 seconds to obtain an approximately equal objective function value of 48.37. 

We observe that in general, the computational time increases as the number of nodes increases. Intuitively, we shall expect that as the number of sprayers increases, the computational time will reduce however, we notice that the computational time increases.

\begin{table}[]
\begin{center}
\caption{Results on the performance of the ALNS metaheuristics.\label{SolverALNS}}
\begin{tabular}{lllllllll}
\multicolumn{1}{c}{} & \multicolumn{3}{c}{Solver} &  & \multicolumn{3}{c}{ALNS}\\ \cline{2-4} \cline{6-8} 
Inst. & Time(s)  & Obj.  & Op. Gap(\%) &    & Time(s)    & Obj.  & Imp(\%) \\ \hline
I-15-S2    & 3888.86  & 48.19 & 0.02    &  & 6.78    & 48.37 & 13.98   \\
I-15-S3    & 1031.30  & 47.91 & 0.00    &  & 9.94    & 48.62 & 14.86   \\
I-17-S2    & 6547.54  & 53.25 & 0.08    &  & 8.28    & 53.14 & 12.60   \\
I-17-S3    & 5030.25  & 55.00 & 0.04    &  & 76.84   & 55.79 & 19.20   \\
I-19-S2    & 7205.01  & 58.42 & 0.13    &  & 138.30  & 58.75 & 14.57   \\
I-19-S3    & 7016.87  & 58.39 & 0.08    &  & 191.54  & 59.73 & 15.13   \\
I-21-S2    & 7206.51  & 63.98 & 0.16    &  & 345.64  & 63.68 & 14.70   \\
I-21-S3    & 7206.92  & 63.99 & 0.14    &  & 782.48  & 64.73 & 20.73   \\
I-23-S2    & 7206.78  & 75.61 & 0.20    &  & 524.47  & 74.54 & 12.54   \\
I-23-S3    & 7206.20  & 70.96 & 0.16    &  & 998.62  & 71.80 & 20.17   \\
I-24-S2    & 7207.73  & 73.30 & 0.21    &  & 541.78  & 72.16 & 13.23   \\
I-24-S3    & 7207.22  & 73.72 & 0.20    &  & 1300.00 & 72.84 & 24.38   \\
I-25-S2    & 7206.57  & 78.93 & 0.22    &  & 226.89  & 77.54 & 10.02   \\
I-25-S3    & 7207.54  & 75.58 & 0.20    &  & 1395.24 & 76.34 & 20.68   \\
I-26-S2    & 7206.88  & 79.47 & 0.23    &  & 857.65  & 78.02 & 14.97   \\
I-26-S3    & 7207.04  & 80.80 & 0.22    &  & 1791.00 & 79.07 & 20.61   \\ \hline
Avg.       & 6424.33  & 66.09 & 0.14    &  & 574.71  & 65.95 & 16.40    \\ 
Min        & 1031.30  & 47.91 & 0.00    &  & 6.78    & 48.37 & 10.02    \\
Max        & 7207.73  & 80.80 & 0.23    &  & 1791.00 & 79.07 & 24.38    \\ \hline
\end{tabular}
\end{center}
\end{table}

\subsubsection{Comparison between models}
In this section, we examine and compare the solutions of the three different models we present in  Section \ref{sec:mathModel} that reflect the preference of the decision maker. For the purposes of illustration, we use instances of size 45 and 60 nodes since these sizes represent a real-life large farm setting across all models with the number of sprayers ranging from 3 to 6. We randomly generate 15 runs for each instance size resulting in a total of 90 experiments. The analysis was on the following measures: total waiting time, total routing time of all sprayers, and the balance of workload across sprayers. The balance of workload across sprayers is measured by comparing the total workload of the busiest sprayer versus the total workload of the least busy sprayer. The total workload of a sprayer is calculated as the routing time, waiting time, spraying (service) time, and refilling time all added up. Then we measure the balance by dividing the workload of the busiest sprayer by the total workload of the least busy sprayer as a measure of the work-load ratio. Such measurement will help in deciding if there is an imbalance in workload across sprayers.

\begin{figure}[!ht]
\centering
\includegraphics[width=14cm]{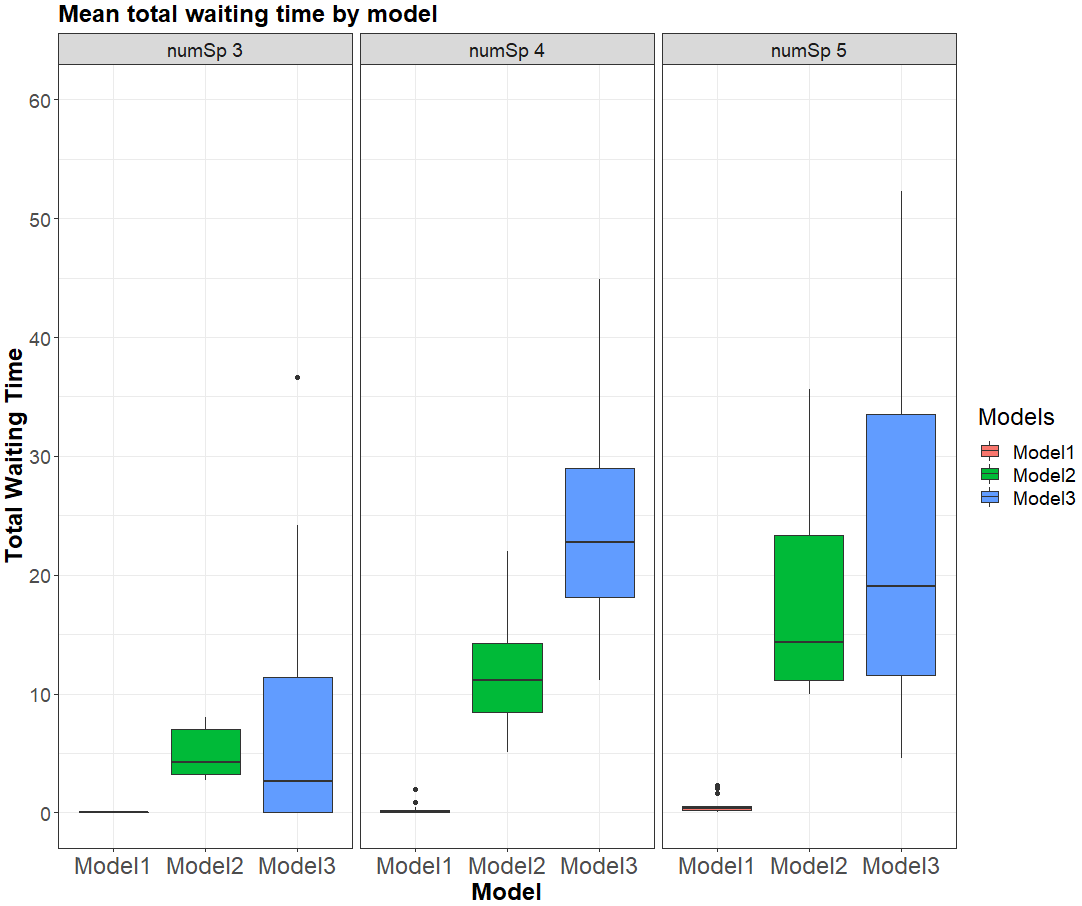}
\caption{Comparison of the total waiting time for different number of sprayers and models}
\label{fig:Waittime}
\end{figure}

Figure \ref{fig:Waittime} illustrates the total waiting time for different number of sprayers $(numSp)$ across the three models. The waiting time of a sprayer at a node is computed as the difference between the arrival time of the tanker and the end of service time of a sprayer at a node. The results show that $Model3$ yields the highest total waiting time followed by $Model2$. $Model1$ yields the lowest waiting time. These results are expected as only $Model1$ focuses on minimizing the waiting time of the sprayers in addition to the routing time. On the other hand, $Model2$ seeks to reduce the latest arrival time of the sprayers while $Model3$ minimizes the total routing time. Overall, we observe an increase in the total waiting time across all the model types as the number of sprayers increases. 

\begin{figure}[!ht]
\centering
\includegraphics[width=14cm]{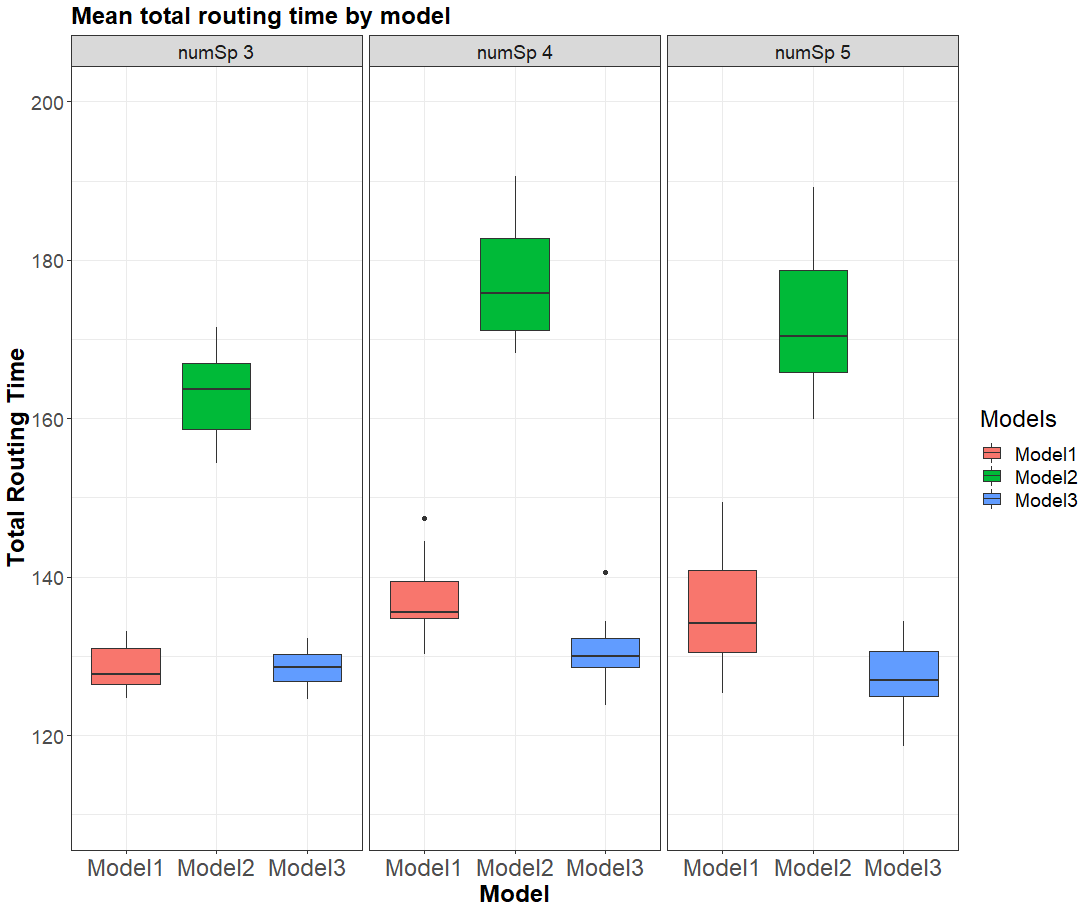}
\caption{Comparison of the total routing time for the models with different number of sprayers}
\label{fig:routingtime}
\end{figure}

Figure \ref{fig:routingtime} depicts the overall routing time of the three models with different number of sprayers $(numSp)$. Obviously, $Model2$ yields the longest total routing time whereas $Model3$ yields the shortest but somewhat lower than $Model1$. The results of $Model3$ and $Model1$ were expected as both models strive to reduce routing time. The total routing time for $Model3$ is somewhat the same as the number of sprayers increases. However, as the number of sprayers increases, we observe an increase in the total routing times of $Model1$ and $Model2$. Without minimizing the routing time we observe that these model approaches attempt to find the optimal routes, whether longest or shortest as long as the routing time does not exceed the $tMax$.

\begin{figure}[!ht]
\centering
\includegraphics[width=14cm]{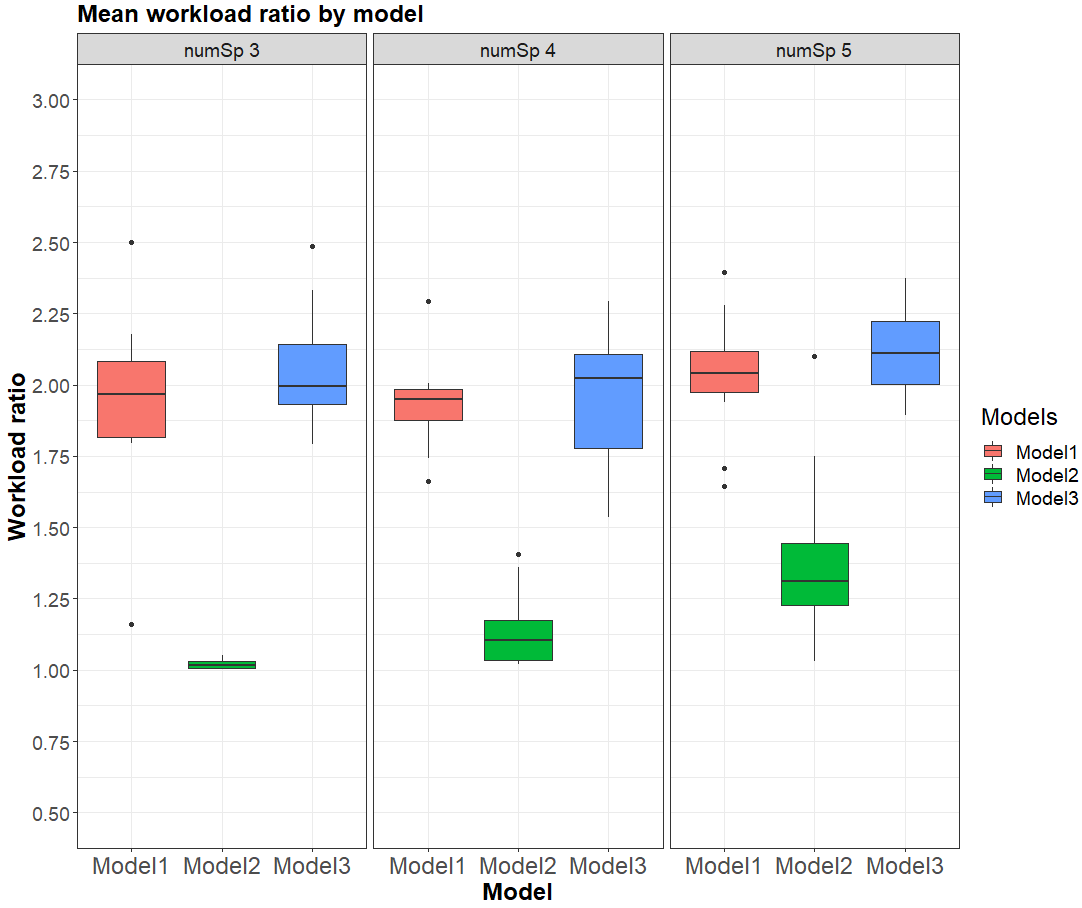}
\caption{Comparison of the workload ratio of the models with different number of sprayers}
\label{fig:workload}
\end{figure}

Lastly, Figure \ref{fig:workload} illustrates the workload ratio under different number of sprayers $numSp$ for the three models. The workload ratio is calculated as the variation between the sprayer with the highest workload and the sprayer with the lowest workload. From the results, we observe that $Model2$ yields the lowest workload ratio. Indicating that there is a great workload balance between the sprayers. For instance, with three sprayers, the workload ratio is approximately 1.00. $Model1$ and $Model3$ yield higher workload with $Model1$ being slightly lower than $Model3$. These results are expected as $Model2$ minimizes the latest arrival of sprayers, thus the model ensures that all sprayers arrive at the depot approximately around the same time. This enables a balanced routing for all the sprayers. We observe that as the number of sprayers increases, the workload ratio increases for $Model2$. But, for $Model1$ and $Model3$ we observe marginal changes in the workload ratio as the number of sprayers increases.\\ 

The benefit of the model comparison study we give extends beyond describing how much each model functions in diverse contexts. Instead, it draws attention to how choosing one model over another causes a significant difference in the performance measures we analyze. For instance, when it comes to measuring the workload ratio and total routing time, the results of $Model1$ and $Model3$ are fairly comparable considering these two performance measures. See Figures \ref{fig:routingtime} and \ref{fig:workload}. However, $Model1$ performs best when the overall waiting time is taken into account Figure\ref{fig:Waittime}. Also, we observe that $Model2$ perform best considering the workload ratio.

\subsubsection{Comparison to Current Practice}
In this section, we compare the performance of the strategy implemented by our partner to solve the sprayer-tanker synchronized routing problem against the optimization approach presented in this study. Our partner company's solution is based on route-first cluster-second strategy presented in Section \ref{sec:Operation}. On the other hand, the optimization approach refers to employing our established metaheuristic to solve a generated instance.

For illustrative purposes, we show the results of comparing the optimization approach to the strategy implemented in practice using randomly generated instances ranging from 30 - 60 nodes. Each setup is run ten times using data regarding the locations of infected areas and the quantity of fertilizer needed at each location and we report the average waiting time and the total routing time under each instance.

The results are shown in Table \ref{Manual}. The format of Table \ref{Manual} is as follows, the first column displays the instance information. The second and third sets of columns show the total waiting time and total routing time under the current practice policy and our optimization approach, respectively. Lastly, the fourth column displays the percentage savings in total time when comparing the optimization approach against the strategy implemented in practice. The savings are calculated as the total time (i.e., routing and waiting times) under the current practice minus the total time under the optimization approach divided by the total time under the current practice multiplied by 100\%. We calculate the mean of mean $Avg.$, minimum  $Min$, and maximum $Max$ of the waiting time and total routing time.

The results show a significant reduction in the waiting time of the sprayers under the optimization approach as opposed to the policy implemented in practice. For example, given instance $I-30-S2$ where the number of nodes is 30 and the number of sprayers is 2, there is no waiting time recorded using the optimization approach but a total waiting time of 3.34 units of time is recorded for solving the same instance using the policy implemented in practice. For the same instance, the total routing time reported for the optimization strategy is 100.83 and that of the current practice is 104.22. We observe that after applying our proposed optimization approach, the overall average percentage time reduction is 5.61\%. Clearly, our optimization approach yields better results than the current practice. We observe that the wait time in both cases increases as the number of sprayers increases. Intuitively, we expect that as the number of sprayers increases, the waiting time should essentially reduce. However, given three sprayers for instance, if the tender tanker is refilling one sprayer, there may be one or two other sprayers waiting for a refill. If it happens that two sprayers are waiting, automatically the third sprayer the tender tanker visits will wait a little longer than the second sprayer causing the overall waiting time to increase. Another interesting finding was that in all cases when the number of sprayers equal three, the overall routing time in our optimization strategy was somewhat longer than the total routing time reported using current practice. However, we notice that in general, our optimization strategy demonstrates excellent performance. These observations emphasize the need to synchronize the tender tanker and sprayer routing so that an optimal solution with the lowest total cost of operation may be reached even without previous knowledge of the sequence and schedule of operations.

Another essential aspect of reducing the overall routing time is minimizing carbon dioxide emissions into the environment. Most companies in the United States want to contribute to a low-carbon economy by lowering their carbon footprint. According to the most recent climate studies, this is critical in order to restrict global warming to 1.5°C, which is necessary to avert the most catastrophic effects of climate change. By implementing our optimization strategy our partner company gets the benefit of reducing carbon emissions.

All of the benefits we discuss above are achieved by our proposed model and algorithm. The management of our partner firm would be quite interested in attaining these advantages by just improving the operational efficiency of their spraying operations without any extra investment.

\begin{table}[]
\begin{center}
\caption{Results of comparing the proposed model against current practice\label{Manual}}
\begin{tabular}{lccccccc}
& \multicolumn{2}{c}{Manual}   & \multicolumn{1}{c}{} & \multicolumn{2}{c}{ALNS}     & \multicolumn{1}{c}{} & \multicolumn{1}{c}{Total } \\ \cline{2-3} \cline{5-6}  
Inst.   & Waiting Time & Total Routing && Waiting Time & Total Routing && Savings (\%) \\ \hline
I-30-S2 & 3.34& 104.22&& 0.00 & 100.83&& 6.26 \\
I-30-S3 & 6.73 & 104.14&& 0.63& 104.54& & 5.14 \\
I-40-S2 & 4.46 & 134.14&& 0.00 & 130.35 & & 5.95\\
I-40-S3 & 10.43& 134.21&& 0.73& 136.85 && 4.89 \\
I-50-S2 & 6.60 & 163.89&& 0.00& 157.12&& 7.84  \\
I-50-S3 & 15.75& 166.63&& 2.00& 173.04 & & 4.02 \\
I-60-S2 & 9.70& 194.99&& 0.00 & 189.03&& 7.65   \\
I-60-S3 & 17.97 & 194.98&& 2.39& 203.86&& 3.15 \\ \hline
Avg.& 9.37& 149.65& & 0.72& 149.45& & 5.61  \\
Min & 3.34& 104.14 & & 0.00& 100.83&& 3.15  \\
Max & 17.97& 194.99&& 2.39& 203.86&& 7.84  \\ \hline
\end{tabular}
\end{center}
\end{table}

\section{Conclusion}
\label{sec:conclude}
In this work, we study a new variant of the vehicle routing problem with multiple synchronization constraints motivated by a real-world farm operating system. We study and model sprayer-tanker synchronized routing problem with multiple trips. We developed three mathematical models based on the key performance criteria that farmers care about (namely, waiting and routing times of sprayers).

The developed mixed integer programming model is an NP-hard and to overcome the computational limitations of solving the model in a reasonable time, we developed a power metaheuristic based on ALNS combined with an intensive local search phase to further improve the quality of the solutions. In comparison to the Gurobi optimization solver, our metaheuristic solves the generated model in a fraction of the time and provides high-quality results. As far as practical value, our optimization approach improves the operational efficiency of the spraying operations by 5.61\%.

An interesting area for future research would be to address stochastic demand because farm managers do not always predict how much fertilizer is needed to spray each infected area with certainty. 

\section*{Compliance with Ethical Standards}
Funding: This study was not funded under any grant.

\section*{Conflict of Interest}

Author Dr. F. A. declares that they have no conflict of interest.

\section*{Ethical approval}
This article does not contain any studies with human participants or animals performed by any of the authors.

\bibliographystyle{pomsref} 

 \let\oldbibliography\thebibliography
 \renewcommand{\thebibliography}[1]{%
 	\oldbibliography{#1}%
 	\baselineskip14pt 
 	\setlength{\itemsep}{10pt}
 }
\bibliography{ref1} 

\begin{thebibliography}{41}
\expandafter\ifx\csname natexlab\endcsname\relax\def\natexlab#1{#1}\fi
\expandafter\ifx\csname url\endcsname\relax
  \def\url#1{{\tt #1}}\fi
\expandafter\ifx\csname urlprefix\endcsname\relax\def\urlprefix{URL }\fi
\expandafter\ifx\csname urlstyle\endcsname\relax
  \expandafter\ifx\csname doi\endcsname\relax
  \def\doi#1{doi:\discretionary{}{}{}#1}\fi \else
  \expandafter\ifx\csname doi\endcsname\relax
  \def\doi{doi:\discretionary{}{}{}\begingroup \urlstyle{rm}\Url}\fi \fi

\bibitem[{Alkaabneh and Diabat(2022)}]{alkaabneh2022multi}
Alkaabneh, Faisal, Ali Diabat. 2022.
\newblock A multi-objective home healthcare delivery model and its solution
  using a branch-and-price algorithm and a two-stage meta-heuristic algorithm.
\newblock {\it Transportation Research Part C: Emerging Technologies\/},
  103838.

\bibitem[{Alkaabneh and Diabat(2023)}]{alkaabneh2023multiobjective}
Alkaabneh, Faisal, Ali Diabat. 2023.
\newblock A multi-objective home healthcare delivery model and its solution
  using a branch-and-price algorithm and a two-stage meta-heuristic algorithm.
\newblock {\it Trnasportation Research Part C: Emerging Technologies\/}, { 147}
  (3), 1-37.

\bibitem[{Alkaabneh et~al.(2023)Alkaabneh, Shehadeh, and
  Diabat}]{alkaabneh2023routing}
Alkaabneh, Faisal, Karmel~S Shehadeh, Ali Diabat. 2023.
\newblock Routing and resource allocation in non-profit settings with equity
  and efficiency measures under demand uncertainty.
\newblock {\it Transportation Research Part C: Emerging Technologies\/}, { 149}
  104023.

\bibitem[{Archetti et~al.(2017)Archetti, Boland, and
  Grazia~Speranza}]{archetti2017matheuristic}
Archetti, Claudia, Natashia Boland, M~Grazia~Speranza. 2017.
\newblock A matheuristic for the multivehicle inventory routing problem.
\newblock {\it INFORMS Journal on Computing\/}, { 29} (3), 377-387.

\bibitem[{Azi et~al.(2014)Azi, Gendreau, and Potvin}]{azi2014adaptive}
Azi, Nabila, Michel Gendreau, Jean-Yves Potvin. 2014.
\newblock An adaptive large neighborhood search for a vehicle routing problem
  with multiple routes.
\newblock {\it Computers \& Operations Research\/}, { 41} 167-173.

\bibitem[{Boysen et~al.(2021)Boysen, Fedtke, and Schwerdfeger}]{boysen2021last}
Boysen, Nils, Stefan Fedtke, Stefan Schwerdfeger. 2021.
\newblock Last-mile delivery concepts: a survey from an operational research
  perspective.
\newblock {\it Or Spectrum\/}, { 43} 1-58.

\bibitem[{Boysen et~al.(2018)Boysen, Schwerdfeger, and
  Weidinger}]{boysen2018scheduling}
Boysen, Nils, Stefan Schwerdfeger, Felix Weidinger. 2018.
\newblock Scheduling last-mile deliveries with truck-based autonomous robots.
\newblock {\it European Journal of Operational Research\/}, { 271} (3),
  1085-1099.

\bibitem[{Cattaruzza et~al.(2016)Cattaruzza, Absi, and
  Feillet}]{cattaruzza2016multi}
Cattaruzza, Diego, Nabil Absi, Dominique Feillet. 2016.
\newblock The multi-trip vehicle routing problem with time windows and release
  dates.
\newblock {\it Transportation Science\/}, { 50} (2), 676-693.

\bibitem[{Cattaruzza et~al.(2018)Cattaruzza, Absi, and
  Feillet}]{cattaruzza2018vehicle}
Cattaruzza, Diego, Nabil Absi, Dominique Feillet. 2018.
\newblock Vehicle routing problems with multiple trips.
\newblock {\it Annals of Operations Research\/}, { 271} (1), 127-159.

\bibitem[{Coindreau et~al.(2021)Coindreau, Gallay, and
  Zufferey}]{coindreau2021parcel}
Coindreau, Marc-Antoine, Olivier Gallay, Nicolas Zufferey. 2021.
\newblock Parcel delivery cost minimization with time window constraints using
  trucks and drones.
\newblock {\it Networks\/}, { 78} (4), 400-420.

\bibitem[{Demir et~al.(2012)Demir, Bekta{\c{s}}, and
  Laporte}]{demir2012adaptive}
Demir, Emrah, Tolga Bekta{\c{s}}, Gilbert Laporte. 2012.
\newblock An adaptive large neighborhood search heuristic for the
  pollution-routing problem.
\newblock {\it European journal of operational research\/}, { 223} (2),
  346-359.

\bibitem[{Drexl(2012)}]{drexl2012synchronization}
Drexl, Michael. 2012.
\newblock Synchronization in vehicle routing—a survey of vrps with multiple
  synchronization constraints.
\newblock {\it Transportation Science\/}, { 46} (3), 297-316.

\bibitem[{Fedtke and Boysen(2017)}]{fedtke2017gantry}
Fedtke, Stefan, Nils Boysen. 2017.
\newblock Gantry crane and shuttle car scheduling in modern rail--rail
  transshipment yards.
\newblock {\it OR spectrum\/}, { 39} (2), 473-503.

\bibitem[{Fran{\c{c}}ois et~al.(2019)Fran{\c{c}}ois, Arda, and
  Crama}]{franccois2019adaptive}
Fran{\c{c}}ois, V{\'e}ronique, Yasemin Arda, Yves Crama. 2019.
\newblock Adaptive large neighborhood search for multitrip vehicle routing with
  time windows.
\newblock {\it Transportation Science\/}, { 53} (6), 1706-1730.

\bibitem[{Fran{\c{c}}ois et~al.(2016)Fran{\c{c}}ois, Arda, Crama, and
  Laporte}]{franccois2016large}
Fran{\c{c}}ois, V{\'e}ronique, Yasemin Arda, Yves Crama, Gilbert Laporte. 2016.
\newblock Large neighborhood search for multi-trip vehicle routing.
\newblock {\it European Journal of Operational Research\/}, { 255} (2),
  422-441.

\bibitem[{Ghilas et~al.(2016)Ghilas, Demir, and
  Van~Woensel}]{ghilas2016adaptive}
Ghilas, Veaceslav, Emrah Demir, Tom Van~Woensel. 2016.
\newblock An adaptive large neighborhood search heuristic for the pickup and
  delivery problem with time windows and scheduled lines.
\newblock {\it Computers \& Operations Research\/}, { 72} 12-30.

\bibitem[{Grangier et~al.(2016)Grangier, Gendreau, Lehu{\'e}d{\'e}, and
  Rousseau}]{grangier2016adaptive}
Grangier, Philippe, Michel Gendreau, Fabien Lehu{\'e}d{\'e}, Louis-Martin
  Rousseau. 2016.
\newblock An adaptive large neighborhood search for the two-echelon
  multiple-trip vehicle routing problem with satellite synchronization.
\newblock {\it European journal of operational research\/}, { 254} (1), 80-91.

\bibitem[{Hashemi~Doulabi et~al.(2020)Hashemi~Doulabi, Pesant, and
  Rousseau}]{hashemi2020vehicle}
Hashemi~Doulabi, Hossein, Gilles Pesant, Louis-Martin Rousseau. 2020.
\newblock Vehicle routing problems with synchronized visits and stochastic
  travel and service times: Applications in healthcare.
\newblock {\it Transportation Science\/}, { 54} (4), 1053-1072.

\bibitem[{Krebs et~al.(2021)Krebs, Ehmke, and Koch}]{krebs2021advanced}
Krebs, Corinna, Jan~Fabian Ehmke, Henriette Koch. 2021.
\newblock Advanced loading constraints for 3d vehicle routing problems.
\newblock {\it OR Spectrum\/}, { 43} (4), 835-875.

\bibitem[{Li et~al.(2020)Li, Wang, Chen, and Bai}]{li2020two}
Li, Hongqi, Haotian Wang, Jun Chen, Ming Bai. 2020.
\newblock Two-echelon vehicle routing problem with time windows and mobile
  satellites.
\newblock {\it Transportation Research Part B: Methodological\/}, { 138}
  179-201.

\bibitem[{Li et~al.(2016)Li, Chen, and Prins}]{li2016adaptive}
Li, Yuan, Haoxun Chen, Christian Prins. 2016.
\newblock Adaptive large neighborhood search for the pickup and delivery
  problem with time windows, profits, and reserved requests.
\newblock {\it European Journal of Operational Research\/}, { 252} (1), 27-38.

\bibitem[{Liu et~al.(2019)Liu, Tao, and Xie}]{liu2019adaptive}
Liu, Ran, Yangyi Tao, Xiaolei Xie. 2019.
\newblock An adaptive large neighborhood search heuristic for the vehicle
  routing problem with time windows and synchronized visits.
\newblock {\it Computers \& Operations Research\/}, { 101} 250-262.

\bibitem[{Mankowska et~al.(2014)Mankowska, Meisel, and
  Bierwirth}]{mankowska2014home}
Mankowska, Dorota~Slawa, Frank Meisel, Christian Bierwirth. 2014.
\newblock The home health care routing and scheduling problem with
  interdependent services.
\newblock {\it Health care management science\/}, { 17} (1), 15-30.

\bibitem[{Meisel and Kopfer(2014)}]{meisel2014synchronized}
Meisel, Frank, Herbert Kopfer. 2014.
\newblock Synchronized routing of active and passive means of transport.
\newblock {\it OR spectrum\/}, { 36} 297-322.

\bibitem[{Mingozzi et~al.(2013)Mingozzi, Roberti, and Toth}]{mingozzi2013exact}
Mingozzi, Aristide, Roberto Roberti, Paolo Toth. 2013.
\newblock An exact algorithm for the multitrip vehicle routing problem.
\newblock {\it INFORMS Journal on Computing\/}, { 25} (2), 193-207.

\bibitem[{Momeni et~al.(2023)Momeni, Mirzapour Al-e Hashem, and
  Heidari}]{momeni2023new}
Momeni, Maryam, SMJ Mirzapour Al-e Hashem, Ali Heidari. 2023.
\newblock A new truck-drone routing problem for parcel delivery by considering
  energy consumption and altitude.
\newblock {\it Annals of Operations Research\/},  1-47.

\bibitem[{Murray and Chu(2015)}]{murray2015flying}
Murray, Chase~C, Amanda~G Chu. 2015.
\newblock The flying sidekick traveling salesman problem: Optimization of
  drone-assisted parcel delivery.
\newblock {\it Transportation Research Part C: Emerging Technologies\/}, { 54}
  86-109.

\bibitem[{Najy et~al.(2023)Najy, Archetti, and Diabat}]{najy2023collaborative}
Najy, Waleed, Claudia Archetti, Ali Diabat. 2023.
\newblock Collaborative truck-and-drone delivery for inventory-routing
  problems.
\newblock {\it Transportation Research Part C: Emerging Technologies\/}, { 146}
  103791.

\bibitem[{Ostermeier et~al.(2023)Ostermeier, Heimfarth, and
  H{\"u}bner}]{ostermeier2023multi}
Ostermeier, Manuel, Andreas Heimfarth, Alexander H{\"u}bner. 2023.
\newblock The multi-vehicle truck-and-robot routing problem for last-mile
  delivery.
\newblock {\it European Journal of Operational Research\/}, { 310} (2),
  680-697.

\bibitem[{Pfeiffer and Schulz(2022)}]{pfeiffer2022alns}
Pfeiffer, Christian, Arne Schulz. 2022.
\newblock An alns algorithm for the static dial-a-ride problem with ride and
  waiting time minimization.
\newblock {\it Or Spectrum\/}, { 44} (1), 87-119.

\bibitem[{Pisinger and Ropke(2007)}]{pisinger2007general}
Pisinger, David, Stefan Ropke. 2007.
\newblock A general heuristic for vehicle routing problems.
\newblock {\it Computers \& operations research\/}, { 34} (8), 2403-2435.

\bibitem[{Quttineh et~al.(2013)Quttineh, Larsson, Lundberg, and
  Holmberg}]{quttineh2013military}
Quttineh, Nils-Hassan, Torbj{\"o}rn Larsson, Kristian Lundberg, Kaj Holmberg.
  2013.
\newblock Military aircraft mission planning: a generalized vehicle routing
  model with synchronization and precedence.
\newblock {\it EURO Journal on Transportation and Logistics\/}, { 2} (1-2),
  109-127.

\bibitem[{Ropke and Pisinger(2006)}]{ropke2006adaptive}
Ropke, Stefan, David Pisinger. 2006.
\newblock An adaptive large neighborhood search heuristic for the pickup and
  delivery problem with time windows.
\newblock {\it Transportation science\/}, { 40} (4), 455-472.

\bibitem[{Sacramento et~al.(2019)Sacramento, Pisinger, and
  Ropke}]{sacramento2019adaptive}
Sacramento, David, David Pisinger, Stefan Ropke. 2019.
\newblock An adaptive large neighborhood search metaheuristic for the vehicle
  routing problem with drones.
\newblock {\it Transportation Research Part C: Emerging Technologies\/}, { 102}
  289-315.

\bibitem[{Schermer et~al.(2019{\natexlab{a}})Schermer, Moeini, and
  Wendt}]{schermer2019hybrid}
Schermer, Daniel, Mahdi Moeini, Oliver Wendt. 2019{\natexlab{a}}.
\newblock A hybrid vns/tabu search algorithm for solving the vehicle routing
  problem with drones and en route operations.
\newblock {\it Computers \& Operations Research\/}, { 109} 134-158.

\bibitem[{Schermer et~al.(2019{\natexlab{b}})Schermer, Moeini, and
  Wendt}]{schermer2019matheuristic}
Schermer, Daniel, Mahdi Moeini, Oliver Wendt. 2019{\natexlab{b}}.
\newblock A matheuristic for the vehicle routing problem with drones and its
  variants.
\newblock {\it Transportation Research Part C: Emerging Technologies\/}, { 106}
  166-204.

\bibitem[{Shaw(1998)}]{shaw1998using}
Shaw, Paul. 1998.
\newblock Using constraint programming and local search methods to solve
  vehicle routing problems.
\newblock {\it International conference on principles and practice of
  constraint programming\/}. Springer, 417-431.

\bibitem[{Soares et~al.(2023)Soares, Marques, Amorim, and
  Parragh}]{soares2023synchronisation}
Soares, Ricardo, Alexandra Marques, Pedro Amorim, Sophie~N Parragh. 2023.
\newblock Synchronisation in vehicle routing: classification schema, modelling
  framework and literature review.
\newblock {\it European Journal of Operational Research\/}, .

\bibitem[{Tamke and Buscher(2021)}]{tamke2021branch}
Tamke, Felix, Udo Buscher. 2021.
\newblock A branch-and-cut algorithm for the vehicle routing problem with
  drones.
\newblock {\it Transportation Research Part B: Methodological\/}, { 144}
  174-203.

\bibitem[{Thomas et~al.(2023)Thomas, Srinivas, and
  Rajendran}]{thomas2023collaborative}
Thomas, Teena, Sharan Srinivas, Chandrasekharan Rajendran. 2023.
\newblock Collaborative truck multi-drone delivery system considering drone
  scheduling and en route operations.
\newblock {\it Annals of Operations Research\/},  1-47.

\bibitem[{Windras~Mara et~al.(2022)Windras~Mara, Norcahyo, Jodiawan,
  Lusiantoro, and Rifai}]{windras2022survey}
Windras~Mara, Setyo~Tri, Rachmadi Norcahyo, Panca Jodiawan, Luluk Lusiantoro,
  Achmad~Pratama Rifai. 2022.
\newblock A survey of adaptive large neighborhood search algorithms and
  applications, .

\end{thebibliography}





\end{document}